\documentclass[12pt]{article}
\usepackage{e-jc}
\usepackage{amsthm,amsmath,amssymb}
\usepackage{amscd}
\usepackage{bm}
\usepackage[breaklinks=true]{hyperref}
\usepackage{tikz}
\usetikzlibrary{shapes}
\tikzstyle{morphism}=[ellipse, fill=green!20!white, draw=black, inner sep=0pt, minimum height=16pt, minimum width=16pt]
\usepackage[backend=bibtex]{biblatex}
\usepackage{graphicx}

\ifpdf
\fi

\newtheorem{theorem}{Theorem}[section]
\newtheorem{defn}[theorem]{Definition}

\newtheorem{cor}[theorem]{Corollary}
\newtheorem{ex}[theorem]{Example}
\newtheorem{prop}[theorem]{Proposition}
\newtheorem{conj}[theorem]{Conjecture}
\newtheorem{lemma}[theorem]{Lemma}
\newtheorem{thm}[theorem]{Theorem}

\newcommand{\bC}{\mathbb{C}}
\newcommand{\bZ}{\mathbb{Z}}
\newcommand{\bN}{\mathbb{N}}
\newcommand{\bQ}{\mathbb{Q}}
\newcommand{\bA}{\mathbb{A}}

\newcommand{\fg}{\mathfrak{g}}

\newcommand{\fsl}{\mathfrak{sl}}
\newcommand{\fgl}{\mathfrak{gl}}

\newcommand{\fS}{\mathfrak{S}}
\newcommand{\brs}{{[\bm\alpha,\bm s]}}

\newcommand{\Hom}{\mathrm{Hom}}
\newcommand{\qbinom}[2]{\genfrac{[}{]}{0pt}{}{#1}{#2}}
\newcommand{\incg}[2][.5in]{\setbox5=\hbox{\;\includegraphics[height=#1]{#2}\;}%
\dimen1=-#1\divide\dimen1 by 2\raise\dimen1\box5}

\newcommand{\bch}{\mathbf{ch}\,}
\newcommand{\fd}{\mathbf{fd}\,}
\newcommand{\wt}{\mathbf{wt}}
\newcommand{\Fix}{\mathrm{Fix}}
\newcommand{\Symm}{\mathrm{Symm}}
\newcommand{\sC}{\mathsf{C}}
\DeclareMathOperator{\id}{id}

\newcommand{\res}{\mathrm{res}}
\newcommand{\bS}{\mathbb{S}}
\newcommand{\tr}{\mathrm{tr}}
\newcommand{\SL}{\mathrm{SL}}
\newcommand{\GL}{\mathrm{GL}}
\newcommand{\Cat}{\mathrm{Cat}}
\newcommand{\cactus}{\mathfrak{C}}	
\title{Invariant tensors and the cyclic sieving phenomenon}
\author{Bruce W. Westbury}
\date{October 2016}

\begin{document}
\maketitle
\begin{abstract}
We construct a large class of examples of the cyclic sieving
phenomenon by exploiting the representation theory of semi-simple Lie algebras.
Let $M$ be a finite dimensional representation of a semi-simple Lie algebra
and let $B$ be the associated Kashiwara crystal. For $r\ge 0$, the triple $(X,c,P)$
which exhibits the cyclic sieving phenomenon is constructed as follows:
the set $X$ is the set of isolated vertices in the crystal $\otimes^rB$;
the map $c\colon X\rightarrow X$ is a generalisation of promotion acting on
standard tableaux of rectangular shape and the polynomial $P$ is the fake degree
of the Frobenius character of a representation of $\mathfrak{S}_r$ related to the natural
action of $\mathfrak{S}_r$ on the subspace of invariant tensors in $\otimes^rM$.
Taking $M$ to be the defining representation of $\mathrm{SL}(n)$ gives the cyclic
sieving phenomenon for rectangular tableaux.
\end{abstract}

\tableofcontents

\section*{Acknowledgement} I would like to thank Martin Rubey for his help and
support in the preparation of this article and the anonymous referee for their time and patience.

I would like to thank the Sage community in general for supporting
this wonderful software and Anne Schilling in particular for implementing crystals
and Travis Scrimshaw for implementing rigged configurations in 
\href{http://www.sagemath.org}{Sage}, \cite{sage}.

\section{Introduction}
In this article we link the combinatorics associated with the cyclic sieving
phenomenon with the representation theory of semisimple Lie algebras and their
associated quantum groups. The cyclic sieving phenomenon is the branch of
algebraic combinatorics which studies pairs $(X,c)$ where $X$ is a finite set
and $c\colon X\rightarrow X$ is a bijection. An invariant of the pair $(X,c)$
is a polynomial $P\in \bZ[q]$ which is only defined modulo $q^r=1$
where $r$ is the order of $c$. This is a complete invariant.

Complex semisimple Lie algebras were classified by Killing in 1888. Their 
representation theory and related theories have been at the forefront of research
ever since. In this article we exploit the representation theory of their
quantised enveloping algebras introduced by Drinfel$'$d and by Jimbo. These quantum
groups were initially introduced to study integrable systems, and two dimensional
vertex models in statistical mechanics in particular. This soon led to applications
to conformal field theory and to low dimensional topology. A subsequent development
consisted in relating the representation theory at roots of unity to the modular
representation theory of algebraic groups. One aspect of this development is the
introduction of the local and global canonical bases by Lusztig and the related
combinatorial theory of crystals. It is this theory that is relevant to this
article.

Each finite dimensional representation has a crystal.
These crystals connect the representation theory with combinatorics.
From the perspective of combinatorics, they extend
the combinatorics associated with partitions,
tableaux and symmetric functions. This combinatorial theory arises
by taking the special case of the Lie algebras $\fgl(n)$ and $\fsl(n)$,
see \cite{MR1338967}.
From the perspective of representation theory, crystals give a combinatorial
interpretation of several branching rules for semisimple Lie algebras.
Prior to the introduction of crystals
these branching rules were expressed as an alternating sum over the Weyl group.
These formula are impractical unless the Weyl group is small.

The first example of this is that
the character of the representation is given by a weighted sum over the
vertices of the crystal. This is a generalisation of the combinatorial
expression of the Schur polynomial as
\begin{equation*}
 s_\lambda (x_1,\dotsc ,x_n) = \sum_T x^{\mathrm{wt}(T)}
\end{equation*}
where the sum is over semistandard tableaux of shape $\lambda$ with entries in
the ordered alphabet $\{1,2,\dotsc ,n\}$.

The most significant applications of crystals arise from the tensor product.
This is a combinatorial rule which gives the crystal of the tensor product
of two representations directly from the crystals of the two representations.
This gives an interpretation of a tensor product multiplicity as the cardinality
of a finite set. This includes, as a special case, the combinatorial interpretation of
Littlewood-Richardson coefficients.

Iterating the tensor product gives the tensor powers of a crystal.
Taking the special case of the vector representation of $\fsl(n)$ gives,
essentially, the Robinson-Schensted correspondence interpreted as a 
combinatorial version of Schur-Weyl duality.

The inspiration for this article was an example of the cyclic sieving phenomenon
first proved in \cite{rhoades}. The set $X$ is the set of standard rectangular
tableaux with $n$ rows and $k$ columns, the bijection $c$ is a combinatorially
defined operation called promotion and the polynomial $P$ is the fake degree of
the Schur function $s_{n^k}$. An alternative proof using webs in the cases $n=2$
and $n=3$ was given in \cite{MR2519848}. A natural question is then whether this
result can be translated into the language of crystals using the above dictionary.
Let $M$ be a representation and $B$ the associated crystal. Then we have $N(\otimes^rM)$,
the subspace of invariant tensors and its combinatorial analogue, $(\otimes^rB)_\ast$,
the set of isolated vertices in $\otimes^rB$. The space $N(\otimes^rM)$ has a natural
action of the symmetric group and so we can define a polynomial by taking the
fake degree of the Frobenius characteristic.

The generalisation of promotion from standard rectangular tableaux to the set
$(\otimes^rB)_\ast$ is new. This is defined as follows,
for $w\in (\otimes^rB)_\ast$,
\begin{itemize}
 \item Remove the first letter of $w$, which is the highest weight of $B$. This
leaves a word $w'$ which is a lowest weight word.
 \item Apply the raising operators to $w'$ to get the highest weight word, $w''$, in the
same component as $w'$.
 \item Add the lowest weight of $B$ to the end of $w''$.
\end{itemize}

An essential ingredient in the proof is a basis of $N(\otimes^rM)$ which is invariant
under the action of the long cycle. In \cite{rhoades} this basis is the Kazhdan-Lusztig
basis and in \cite{MR2519848} this is the web basis introduced in \cite{MR1403861}.
The web basis is simpler but is
only known in a small number of examples. For the general case we use the basis
constructed in \cite[\S 27.3]{MR1227098} using the theory of based modules. This agrees
with the Kazhdan-Lusztig basis for rectangular tableaux.

These definitions allow us to formulate a generalisation of the motivating example
of the cyclic sieving phenomenon. This does not involve any explicit mention of the
quantised enveloping algebras. However the proof does make essential use of their
representation theory. The problem is to relate two actions of the cyclic group;
one is the action of the long cycle on $N(\otimes^rM)$ which is defined using the
natural action of the symmetric group; the other is the promotion operator acting
on $(\otimes^rB)_\ast$ and which is defined combinatorially.

The quantised enveloping algebra is used to relate these two cyclic group actions.
This introduces a parameter $q$ and, loosely speaking, the action of the long cycle
is defined for $q=1$ and the crystals are defined for $q^{-1}=0$; and so we are
interpolating between these. For $q=1$ the symmetric group acts on tensor powers
by permuting indices. This no longer holds for quantum groups and Drinfel$'$d
introduced two weakenings. The more familiar weakening replaces the symmetric
groups by the braid groups; this only plays a minor role in this paper. The other
weakening replaces the symmetric groups by the cactus groups. This weakening is
crucial to this paper because it passes to crystals (whereas the braiding does not).
For more background on these structures see \cite{MR2497964}.

Section \ref{sec:g2p} gives some examples for the seven dimensional representation of $G_2$.
This is included to illustrate the theory that has been developed in this paper.

In this article we make extensive use of the string diagram notation for
tensors. This notation was introduced in \cite{MR776784} and used
throughout \cite{MR2418111}. This notation became widely accepted when
braided monoidal categories became popular and has antecedents in
the string diagrams for the permutation groups and braid groups,
in \cite{MR1503378} and in \cite{MR0127765}.

Our notation for symmetric functions is standard and follows \cite{MR1354144}
and \cite{MR1676282}. For a partition $\lambda$, $h_\lambda$ is the complete
homogeneous symmetric function, $e_\lambda$ is the elementary symmetric function,
$p_\lambda$ is the power sum symmetric function and $s_\lambda$ is the Schur function.
We denote the space of homogeneous symmetric functions of degree $r$ by $\Symm(r)$.
The algebra of symmetric functions has an involution $\omega$ which is determined by
$\omega(h_\lambda)=e_\lambda$, $\omega(e_\lambda)=h_\lambda$ and
$\omega(s_\lambda)=s_\lambda'$ where $\lambda'$ is the partition conjugate to $\lambda$.

The $q$-integers, $q$-factorials and $q$-binomial coefficients are used in \S \ref{sec:csp}
follow the standard combinatorial conventions.
These are elements of $\bZ[q]$, so are polynomials. The $q$-integer $[r]$ is 
\begin{equation*}
 [r] = \frac{1-q^r}{1-q}
\end{equation*}
and then $[r]!=[r][r-1]\dotsb [1]$ and $\qbinom{r}{s}=\frac{[r]!}{[s]![r-s]!}$.
There is a potential source of confusion. Although no $q$-integers
appear in \S \ref{sec:qg}, in the literature on quantised enveloping algebras, the $q$-integer
$[r]$ is
\begin{equation*}
 [r] = \frac{q^r-q^{-r}}{q-q^{-1}}
\end{equation*}
The $q$-factorials and $q$-binomial coefficients are then defined in the same way.
These are elements of $\bZ[q,q^{-1}]$ and so Laurent polynomials. These are all
invariant under the involution $q\leftrightarrow q^{-1}$.

Examples of the cyclic sieving phenomenon which are obtained using
Theorem \ref{thm:fund} are given in \cite{3592}. These examples arise
in the representation theory of the symplectic groups.
\section{Cyclic sieving phenomenon}\label{sec:csp}
The \emph{cyclic sieving phenomenon} was introduced in \cite{RSW}
and there are expositions in \cite{MR2866734} and \cite{MR3156682}.
These well-written articles emphasise the combinatorial aspects. Here we review
the cyclic sieving phenomenon from the perspective of representation
theory.

\subsection{Cyclic groups}
Put $\omega=\exp(2\pi i/r)$. Let $C_r$ be a cyclic group of order $r$
with a generator $c$. Then the irreducible complex representations
of $C_r$ are one dimensional and are given by $c^p \mapsto \omega^{pk}$
for $0\le k < r$. This gives an identification of the representation
ring of $C_r$ with $\bZ[q]/\left\langle q^r-1 \right\rangle$. The element
$P\in \bZ[q]/\left\langle q^r-1\right\rangle$ associated to a representation is
characterised by the property that, for all $k$,
\begin{equation*}
 \tr(c^k)=P(\omega^k)
\end{equation*}

A permutation representation of $C_r$ is a set $X$ with a bijection
$c\colon X\rightarrow X$ of order $r$. This gives a linear representation
and hence an element $P\in\bZ[q]/\left\langle q^r-1 \right\rangle$.
For example, if $d|r$ then there is a transitive permutation representation
of size $r/d$ such that the stabiliser of any point is the subgroup of order
$d$ generated by $c^{r/d}$. The associated polynomial is $(q^r-1)/(q^d-1)$.
This gives all the transitive permutation representations.
It follows that the polynomial
associated to a permutation representation determines the permutation representation.
A more explicit formulation is that the coefficient of $q^k$ is the number of
orbits such that the order of the stabiliser of any point in the orbit divides $k$.

\begin{ex}\label{pm}
For each $r$ we have the set, $X(r)$, of noncrossing perfect matchings on the set
$[2r]=\{1,2,\dotsc ,2r\}$. The cardinality of $X(r)$ is the Catalan number
\begin{equation*}
 \Cat(r) = \frac1{(r+1)}\binom{2r}{r}
\end{equation*}
The diagram of a noncrossing perfect matching consists of $r$ arcs drawn
in the unit disk in the plane. The element  $u\in [2r]$ is identified with
the boundary point $(\cos(u\pi/r),\sin(u\pi/r))$. For each pair $(u,v)$ in the
perfect matching there is an arc connecting the boundary points $u$ and $v$.
These arcs are noncrossing.

Rotation by the angle $\pi/r$ defines a bijection $c\colon X(r)\rightarrow X(r)$
of order $2r$.

Define a $q$-analogue of the Catalan number by
\begin{equation*}
 \Cat_q(r) = q^{r(r-1)}\frac1{[r+1]}\qbinom{2r}{r}
\end{equation*}

Then $(X(r),c,\Cat_q(r))$ exhibits the cyclic sieving phenomenon.

These polynomials and their reductions $\mod q^{2r}-1$ are:
\begin{center}\begin{tabular}{r|r|r}
0 & 1 & 1 \\
1 & 1 & 1 \\
2 & $q^4+q^2$ & $q^2+1$ \\
3 & $q^{12}+q^{10}+q^9+q^8+q^6$ & $q^4+q^3+q^2+2$
\end{tabular}\end{center}
For $r=2$ this says that there is one orbit of order two. For $r=3$ this says that there
is one orbit of order two and one of order three.

The orbit of order two is shown in Figure \ref{fig:o2}.
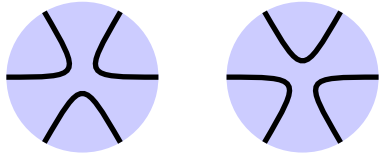
\begin{figure}
\begin{center}
\begin{tikzpicture}[line width=2pt]
\fill[color=blue!20] (0,0) circle (1);
\draw (0:1) .. controls (0,0) and (0,0) .. (60:1);
\draw (120:1) .. controls (0,0) and (0,0) .. (180:1);
\draw (240:1) .. controls (0,0) and (0,0) .. (300:1);
\end{tikzpicture}\qquad 
\begin{tikzpicture}[line width=2pt]
\fill[color=blue!20] (0,0) circle (1);
\draw (60:1) .. controls (0,0) and (0,0) .. (120:1);
\draw (180:1) .. controls (0,0) and (0,0) .. (240:1);
\draw (300:1) .. controls (0,0) and (0,0) .. (360:1);
\end{tikzpicture}
\end{center}
\caption{Orbit of order two}\label{fig:o2}
\end{figure}
The orbit of order three is shown in Figure \ref{fig:o3}.

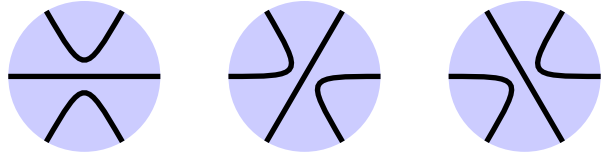
\begin{figure}
\begin{center}
\begin{tikzpicture}[line width=2pt]
\fill[color=blue!20] (0,0) circle (1);
\draw (0:1) .. controls (0,0) and (0,0) .. (180:1);
\draw (60:1) .. controls (0,0) and (0,0) .. (120:1);
\draw (240:1) .. controls (0,0) and (0,0) .. (300:1);
\end{tikzpicture}\qquad 
\begin{tikzpicture}[line width=2pt]
\fill[color=blue!20] (0,0) circle (1);
\draw (60:1) .. controls (0,0) and (0,0) .. (240:1);
\draw (120:1) .. controls (0,0) and (0,0) .. (180:1);
\draw (300:1) .. controls (0,0) and (0,0) .. (360:1);
\end{tikzpicture}\qquad 
\begin{tikzpicture}[line width=2pt]
\fill[color=blue!20] (0,0) circle (1);
\draw (120:1) .. controls (0,0) and (0,0) .. (300:1);
\draw (180:1) .. controls (0,0) and (0,0) .. (240:1);
\draw (360:1) .. controls (0,0) and (0,0) .. (60:1);
\end{tikzpicture}
\end{center}
\caption{Orbit of order three}\label{fig:o3}
\end{figure}
\end{ex}

The following is a basic example of the cyclic sieving phenomenon.
This is \cite[Theorem 1.6]{RSW} and there are several proofs given.
The special case in which $\lambda$ has two parts is \cite[Theorem 2.1]{MR2866734}.
\begin{lemma}\label{lem:csp} Let $\lambda\vdash r$. Let $X_\lambda$ be the set of set decompositions 
of $\{1,2,\dotsc ,r\}$ into blocks of sizes $\lambda_1,\lambda_2,\dotsc$.
This is a permutation representation of $\fS_r$ and we take $c$ to be the action
of the long cycle. Let $P_\lambda$ be the $q$-multinomial coefficient
\begin{equation*}
 \qbinom{r}{\lambda} = \frac{[r]!}{\prod_i[\lambda_i]! }
\end{equation*}
Then $(X_\lambda,c,P_\lambda)$ exhibits the cyclic sieving phenomenon.
\end{lemma}

\subsection{Symmetric groups}
Next we explain how the representation theory of the symmetric group can be used
to compute the polynomial of a representation of $C_r$.

Let $\Gamma$ be a finite group. Let $K(\Gamma)$ be the representation ring of
complex representations of $\Gamma$; so each representation $V$ gives an element $[V]\in K(\Gamma)$.
Let $c\in \Gamma$ be an element of order $r$ defined
up to conjugacy. Then the restriction of representations induces a map
$\res\colon K(\Gamma)\rightarrow \bZ[q]/\left\langle q^r-1 \right\rangle$.

\begin{lemma}
Let $V$ be a representation of $\Gamma$ with a basis $X$ which is fixed by $c$.
Put $P=\res([V])$. Then $P$ is the character of $(X,c)$.
\end{lemma}

In this article we take $\Gamma$ to be the symmetric group, $\fS_r$, and $c\in\fS_r$
to be the long cycle. First we identify $K(\fS_r)$ with homogeneous symmetric
functions of degree $r$ using the characteristic map. The \emph{Frobenius character}
 of a linear representation, $V$, is
given in \cite[\S I.7 (7.2)]{MR1354144} and \cite[\S 7.18]{MR1676282}.
\begin{equation*}
 \bch(V) = \frac{1}{r!} \sum_{w\in \fS_r} \tr(w) p_{\lambda(w)}
\end{equation*}
where $\lambda(w)$ is the cycle type of $w$ and $\tr(w)$ is the matrix trace of $w$ acting on $V$.

The \emph{cycle index series} of a permutation representation was introduced
in \cite{MR1506633} and is the main topic of \cite{BLL}. The cycle index series
of a permutation representation, $X$, is
\begin{equation*}
 Z(X) = \frac{1}{r!} \sum_{w\in \fS_r} \Fix(w) p_{\lambda(w)}
\end{equation*}
where $\Fix(w)$ is the number of points of $X$ fixed by $w$.

Let $X$ be a permutation representation and $\bC X$ the associated linear representation.
Then the equality $\bch(\bC X)=Z(X)$ follows from the simple observation that the
matrix trace of $w$ acting on $\bC X$ is the number of points of $X$ fixed by the action
of $w$.

The \emph{fake degree} is a linear map from symmetric functions to $\bZ[q]$.
Let $f$ be a symmetric function. Then the stable principal specialisation
of $f$ is the formal power series $\mathbf{ps}(f)=f(1,q,q^2,\dotsc )$.
If $f$ is homogeneous of degree $r$ then the fake degree of $f$ is
\begin{equation}\label{fd:ps}
 \fd(f) = \frac{\mathbf{ps}(f)}{[r]!}
\end{equation}

Let $\lambda$ be a partition and put $r=|\lambda|$.
The polynomial $\fd(h_\lambda)$ is given in \cite[7.8.3 Proposition]{MR1676282}. We have that
$\fd(h_\lambda)$ is the $q$-multinomial coefficient
\begin{equation*}
 \fd(h_\lambda) = \qbinom{r}{\lambda}=\frac{[r]!}{\prod_i [\lambda_i]!}
\end{equation*}

There are two expressions for $\fd(s_\lambda)$.
These are given in Corollaries 7.21.3 and 7.21.5 in
\cite{MR1676282}. One is a $q$-analogue of the hook length
formula and is
\begin{equation}\label{hook}
 \fd(s_\lambda) = q^{b(\lambda)}\frac{[r]!}{\prod_{(i,j)\in\lambda} [h(i,j)]}
\end{equation}
where $(i,j)$ is a cell of $\lambda$, $h(i,j)$ is its hook length
and $b(\lambda)=\lambda_2+2\lambda_3+3\lambda_4+\dotsb$. The other is
\begin{equation}\label{fd:maj}
 \fd(s_\lambda) = \sum_T q^{\mathrm{maj}(T)}
\end{equation}
where the sum is over standard tableaux of shape $\lambda$ and $\mathrm{maj}(T)$
is the major index of the tableau $T$.

The result that \eqref{fd:ps} and \eqref{fd:maj} agree is given in
\cite[Proposition 2.1]{MR732346}. One way to see this is to observe
that both definitions have the property that if $f\in\Symm(r)$
and $g\in\Symm(s)$ then
\begin{equation*}
 \fd(fg) = \qbinom{r+s}{r} \fd(f) \fd(g)
\end{equation*}
This uses the extension of \eqref{fd:maj} to skew shapes.

\begin{prop}
Put $\Gamma=\fS_r$ and take $c\in\fS_r$ to be the long cycle.
Identify $K(\fS_r)$ with $\Symm(r)$ using the Frobenius character so we
have the restriction map
$\res\colon \Symm(r) \rightarrow \bZ[q]/\left\langle q^r-1 \right\rangle$.
Then $\res(f) = \fd(f) \mod q^r-1$ for all $f\in\Symm(r)$.
\end{prop}

\begin{proof}[First proof]
Since $Symm(r)$ is a free abelian group, and both maps are linear,
it is sufficient to check the result on a basis.
Take the complete homogeneous symmetric functions
as basis and use Lemma \ref{lem:csp}.
\end{proof}

\begin{proof}[Second proof]
Since $\Symm(r)$ is a free abelian group, and both maps are linear,
it is sufficient to check the result on the basis of Schur functions.

For $0\le k<r$ and $\lambda\vdash r$, let $T(\lambda,k)$ be the number of
standard tableaux of shape $\lambda$ whose major index equals $k$ mod $r$.

It is shown in \cite[Theorem 8.6]{MR1039352}, \cite{MR1867283}, \cite[Theorem 8.8]{MR1231799} that
\begin{equation*}
[ \mathrm{Ind}_{C_r}^{\fS_r}(q^k) : S(\lambda) ] = T(\lambda,k)
\end{equation*}

Frobenius reciprocity is, essentially, the result that induction
is left adjoint to restriction. This gives
\begin{equation*}
[ \mathrm{Ind}_{C_r}^{\fS_r}(q^k) : S(\lambda) ] =
[ \mathrm{Res}_{C_r}^{\fS_r}(S(\lambda)) : q^k ] 
\end{equation*}

Hence
\begin{equation*}
[ \mathrm{Res}_{C_r}^{\fS_r}(S(\lambda)) : q^k ] = T(\lambda,k)
\end{equation*} 
\end{proof}

A straightforward consequence is:
\begin{thm}\label{thm:fund} Let $V$ be a linear representation of $\fS_r$ with a basis $X$
which is preserved by the long cycle, $c$. Put $P=\fd(\bch(V))$.
Then $(X,c,P)$ exhibits the cyclic sieving phenomenon.
\end{thm}

\begin{cor} Let $X$ be a permutation representation of $\fS_r$.
Put $P=\fd(Z(X))$. Then $(X,c,P)$ exhibits the cyclic sieving phenomenon.
\end{cor}

\begin{ex}
This is a continuation of Example \ref{pm}.
Let $N(r)$ be the vector space with basis $X(r)$. Then we define an action of the
symmetric group $\fS_{2r}$ on $N(r)$. First we draw the diagram of a noncrossing
perfect matching in the upper half plane. Then the action of a permutation
is given by stacking the diagram of the perfect matching on top of the string
diagram of the permutation and then using the relations
\[ \incg{csp.4} = -2 \incg{csp.48} \]
\[ \incg{csp.49} = - \incg{csp.50} - \incg{csp.51} \]

Then this action has the property that rotation is given by the action of the long cycle
in $\fS_{2r}$. This representation of $\fS_{2r}$ is irreducible and the Frobenius character
is the Schur function $s_{2^r}$ (where $2^r$ is the partition conjugate to $[r,r]$).
Putting $\lambda=2^r$ in \eqref{hook} gives $\fd(s_{2^r})=\Cat_q(r)$.
\end{ex}

\section{Monoidal categories}
In this section we give a concise account of monoidal
categories.  Here we only discuss the strict versions of
these concepts; this is for brevity and simplicity.

We start with monoidal categories and then impose several different
additional structures. The most familiar additional structure gives
symmetric monoidal categories. There are two different ways to weaken
this namely, braided monoidal categories and coboundary monoidal
categories. These were both introduced by  in his deformation
quantisation approach to quantum groups. In Jimbo's alternative approach
to quantum groups using finite presentations the braided structure is
implemented using the universal $R$-matrix and the coboundary structure
was described more recently in \cite{MR2219257}.

\subsection{Monoidal categories}
Monoidal and symmetric monoidal categories were introduced in \cite{MR0170925}
and a coherence theorem is given in \cite{MR593254}. A more recent reference
is \cite{MR1268782}. In this paper we use the diagram calculus for
monoidal and symmetric monoidal categories described in \cite[Sections 1,2]{MR1113284}.

\begin{defn}
A \emph{strict monoidal category} consists of the following data:
\begin{itemize}
 \item a category $\sC$,
 \item a functor $\otimes \colon \sC \times \sC\rightarrow\sC$,
 \item an object $I\in\sC$.
\end{itemize}
The functor $\otimes$ is required to be associative. This is the
condition that
\begin{equation*}
\otimes \circ (\otimes \times\id) = \otimes \circ (\id \times \otimes)
\end{equation*}
as functors $\sC \times \sC \times \sC\rightarrow\sC$. Equivalently,
the following diagram commutes:
\begin{equation*}\begin{CD}
\sC \times \sC \times \sC @>{\otimes \times\id}>> \sC \times \sC \\
@V{\id \times \otimes}VV @VV{\otimes}V \\
\sC \times \sC @>>{\otimes}> \sC
\end{CD}\end{equation*}
The conditions on the object $I$ are that the functors $-\otimes I$ and $I\otimes-$
are both equivalences of categories.
\end{defn}

Let $\sC$ and $\sC'$ be strict monoidal categories. A functor
$F\colon\sC\rightarrow\sC'$ is \emph{monoidal}
if $F(I)=I'$ and the following diagram commutes
\begin{equation*}\begin{CD}
\sC \times \sC  @>{F\times F}>> \sC' \times \sC' \\
@V{\otimes}VV @VV{\otimes}V \\
\sC @>>{F}> \sC'
\end{CD}\end{equation*}

In this article we will make extensive use of string diagrams for
morphisms in a monoidal
category. A morphism is represented by a graph embedded in a
rectangle with boundary points
on the top and bottom edges. The edges of the graph are directed and
are labelled by objects.
The vertices are labelled by morphisms. A morphism $f\colon
x\rightarrow y$ is drawn as follows:
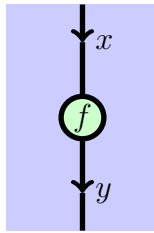
\begin{equation*}
\begin{tikzpicture}[line width=2pt]
\fill[color=blue!20] (-1,-1.5) rectangle (1,1.5);
\draw (0,1.5)[->] -- (0,1) node[anchor=west] {$x$}; \draw[->] (0,1) -- 
(0,-1);
\draw (0,-1) node[anchor=west] {$y$} -- (0,-1.5);
\draw (0,0) node[morphism] {$f$};
\end{tikzpicture}
\end{equation*}

\begin{figure}
\begin{center}
\begin{tikzpicture}[line width=2pt]
\fill[color=blue!20] (-1,-1.75) rectangle (1.2,1.75);
\draw (0,1.75)[->] -- (0,1) node[anchor=west] {$x$};
\draw[->] (0,1) -- (0,-1);
\draw (0,-1) node[anchor=west] {$z$} -- (0,-1.75);
\draw (0,0) node[morphism] {$g\circ f$};
\draw (2,0) node {$=$};
\fill[color=blue!20] (3,-1.75) rectangle (5,1.75);
\draw (4,1.75)[->] -- (4,1.333);
\draw (4,1.333)[->] node[anchor=west] {$x$} -- (4,0);
\draw (4,0)[->] node[anchor=west] {$y$} -- (4,-1.333);
\draw (4,-1.333) node[anchor=west] {$z$} -- (4,-1.75);
\draw (4,0.666) node[morphism] {$f$};
\draw (4,-0.666) node[morphism] {$g$};
\end{tikzpicture}\qquad
\begin{tikzpicture}[line width=2pt]
\fill[color=blue!20] (-1,-1.5) rectangle (1.2,1.5);
\draw (0,1.5)[->] -- (0,1) node[anchor=west] {$x\otimes u$};
\draw[->] (0,1) -- (0,-1);
\draw (0,-1) node[anchor=west] {$y\otimes v$} -- (0,-1.5);
\draw (0,0) node[morphism] {$f\otimes g$};
\draw (2,0) node {$=$};
\fill[color=blue!20] (3,-1.5) rectangle (5,1.5);
\draw (3.65,1.5)[->] -- (3.65,1) node[anchor=east] {$x$};
\draw[->] (3.65,1) -- (3.65,-1);
\draw (3.65,-1) node[anchor=east] {$y$} -- (3.65,-1.5);
\draw (3.65,0) node[morphism] {$f$};
\draw (4.35,1.5)[->] -- (4.35,1) node[anchor=west] {$u$};
\draw[->] (4.35,1) -- (4.35,-1);
\draw (4.35,-1) node[anchor=west] {$v$} -- (4.35,-1.5);
\draw (4.35,0) node[morphism] {$g$};
\end{tikzpicture}
\end{center}
\caption{Composition and tensor product}\label{fig:comp}
\end{figure}
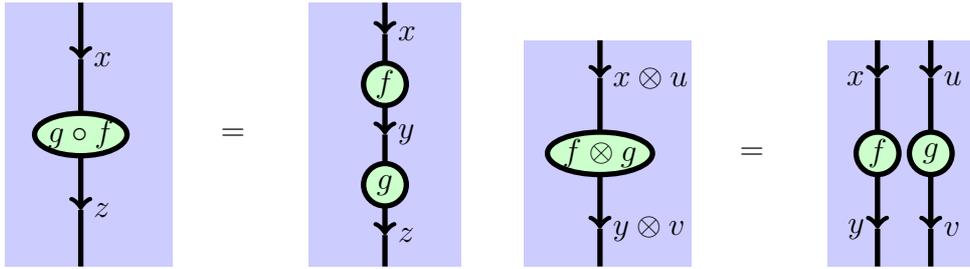

The rules for combining diagrams are that composition is given by
stacking diagrams and the tensor product is given by putting diagrams
side by side.  The corresponding string diagrams are shown in
Figure~\ref{fig:comp}.

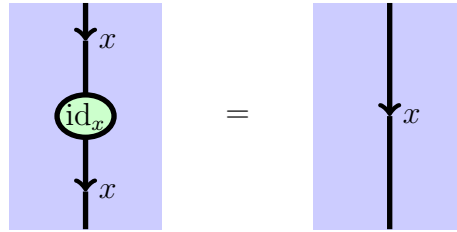
\begin{figure}
\begin{center}
\begin{tikzpicture}[line width=2pt]
\fill[color=blue!20] (-1,-1.5) rectangle (1,1.5);
\draw (0,1.5)[->] -- (0,1) node[anchor=west] {$x$}; \draw[->] (0,1) -- 
(0,-1);
\draw (0,-1) node[anchor=west] {$x$} -- (0,-1.5);
\draw (0,0) node[morphism] {$\id_x$};
\draw (2,0) node {$=$};
\fill[color=blue!20] (3,-1.5) rectangle (5,1.5);
\draw (4,1.5)[->] -- (4,0); \draw (4,0) node[anchor=west] {$x$} -- (4,-1.5);
\end{tikzpicture}
\end{center}
\caption{Identity morphism}\label{fig:id}
\end{figure}

Rules for simplifying diagrams are that vertices labelled by an
identity morphism can be omitted and edges labelled by the trivial
object, $I$, can be omitted.  The string diagram for the rule that an
identity morphism can be omitted is shown in Figure~\ref{fig:id}.

\subsection{Symmetric monoidal categories}
The examples of symmetric monoidal categories which motivated the definition
are categories of representations of groups and of Lie algebras.
\begin{defn} A \emph{strict symmetric monoidal category} is a strict
monoidal category $\sC$ together with
natural isomorphisms $\alpha(x,y)\colon x\otimes y\rightarrow
y\otimes x$ for all objects $x,y$.
These are required to satisfy the conditions that 
\begin{equation}\label{a:square}
\alpha(x,y)\circ\alpha(y,x) = \id_{y\otimes x}
\end{equation}
and that the following diagrams commute for all $f\colon u\rightarrow v$ and all $x,y,z$:
\begin{equation}\label{a:slide1}\begin{CD}
u\otimes x @>{\alpha(u,x)}>> x\otimes u \\
@V{f\otimes \id_x}VV @VV{\id_x\otimes f}V \\
v\otimes x @>>{\alpha(v,x)}> x\otimes v
\end{CD}\end{equation}
\begin{equation}\label{a:slide2}\begin{CD}
x\otimes u @>{\alpha(x,u)}>> u\otimes x \\
@V{\id_x\otimes f}VV @VV{f\otimes \id_x}V \\
x\otimes v @>>{\alpha(x,v)}> v\otimes x
\end{CD}\end{equation}
\begin{equation}\label{a:tensor1}\begin{CD}
x\otimes y\otimes z @>{\alpha(x,y)\otimes\id_z}>> y\otimes x\otimes z\\
@V{\alpha(x,y\otimes z)}VV @VV{\id_y\otimes \alpha(x,z)}V \\
y\otimes z\otimes x @= y\otimes z\otimes x
\end{CD}\end{equation}
\begin{equation}\label{a:tensor2}\begin{CD}
x\otimes y\otimes z @>{\id_x\otimes\alpha(y,z)}>> y\otimes x\otimes z\\
@V{\alpha(x\otimes y, z)}VV @VV{ \alpha(x,z)\otimes\id_y}V \\
z\otimes x\otimes y @= z\otimes x\otimes y
\end{CD}\end{equation}
\end{defn}

A \emph{symmetric monoidal functor} is a functor compatible with this
structure.

The string diagram for $\alpha(x,y)$ is as follows:
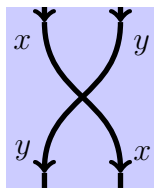
\begin{equation*}
\begin{tikzpicture}[line width=2pt]
\fill[color=blue!20] (-1,-1.2) rectangle (1,1.2);
\draw (-0.5,-1) -- (-0.5,-1.2);
\draw (0.5,-1) -- (0.5,-1.2);
\draw (-0.5,1.2)[->] -- (-0.5,1);
\draw (0.5,1.2)[->] -- (0.5,1);
\draw (-0.5,1)[->] node[anchor=north east] {$x$}  .. controls (-0.5,0) and (0.5,0) .. (0.5,-1) node[anchor=south west] {$x$};
\draw (0.5,1)[->] node[anchor=north west] {$y$}  .. controls (0.5,0) and (-0.5,0) .. (-0.5,-1) node[anchor=south east] {$y$};
\end{tikzpicture}
\end{equation*}
\begin{figure}
\begin{center}
\begin{tikzpicture}[line width=2pt]
\fill[color=blue!20] (-1,-1.2) rectangle (1,3.2);
\draw (-0.5,-1) -- (-0.5,-1.2);
\draw (0.5,-1) -- (0.5,-1.2);
\draw (-0.5,3.2)[->] -- (-0.5,3) node[anchor=north east] {$x$};
\draw (0.5,3.2)[->] -- (0.5,3)  node[anchor=north west] {$y$};
\draw (-0.5,1)[->]   .. controls (-0.5,0) and (0.5,0) .. (0.5,-1)
node[anchor=south west] {$y$};
\draw (0.5,1)[->]  .. controls (0.5,0) and (-0.5,0) .. (-0.5,-1)
node[anchor=south east] {$x$};
\draw (-0.5,3)[->]  .. controls (-0.5,2) and (0.5,2) .. (0.5,1)
node[anchor= west] {$x$};
\draw (0.5,3)[->]  .. controls (0.5,2) and (-0.5,2) .. (-0.5,1)
node[anchor= east] {$y$};
\draw (1.5,1) node {$=$};
\fill[color=blue!20] (2,-1.2) rectangle (4,3.2);
\draw (2.5,3.2)[->] -- (2.5,1);
\draw (2.5,1) node[anchor= east] {$x$} -- (2.5,-1.2);
\draw (3.5,3.2)[->] -- (3.5,1);
\draw (3.5,1) node[anchor= west] {$y$} -- (3.5,-1.2);
\end{tikzpicture}
\end{center}
\caption{Symmetry}\label{fig:sym}
\end{figure}
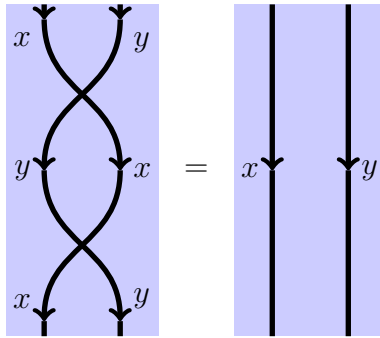
The string diagram for condition~\eqref{a:square} is shown in
Figure~\ref{fig:sym}, the string diagrams for
conditions~\eqref{a:slide1} and~\eqref{a:slide2} are shown in
Figure~\ref{fig:rels-slide}, and those for
conditions~\eqref{a:tensor1} and~\eqref{a:tensor2} are shown in
Figure~\ref{fig:rels-tensor}.
\begin{figure}
\begin{center}
\begin{tikzpicture}[line width=2pt]
\fill[color=blue!20] (-1,-1.2) rectangle (1,3.2);
\draw (-0.5,-1) node[anchor=south east] {$x$} -- (-0.5,-1.2);
\draw (0.5,-1) node[anchor=south west] {$v$} -- (0.5,-1.2);
\draw (-0.5,3.2)[->] -- (-0.5,1);
\draw (0.5,3.2)[->] -- (0.5,1);
\draw (-0.5,3.2)[->] -- (-0.5,3) node[anchor=north east] {$u$};
\draw (0.5,3.2)[->] -- (0.5,3)  node[anchor=north west] {$x$};
\draw (-0.5,1)[->] node[anchor=east] {$v$}  .. controls (-0.5,0) and
(0.5,0) .. (0.5,-1) ;
\draw (0.5,1)[->]  node[anchor=west] {$x$} .. controls (0.5,0) and
(-0.5,0) .. (-0.5,-1);
\draw (-0.5,2) node[morphism] {$f$};
\draw (1.5,1) node {$=$};
\fill[color=blue!20] (2,-1.2) rectangle (4,3.2);
\draw (2.5,3.2)[->] -- (2.5,3) node[anchor=north east] {$u$};
\draw (3.5,3.2)[->] -- (3.5,3)  node[anchor=north west] {$x$};
\draw (2.5,1) node[anchor=  east] {$x$} -- (2.5,-1.2);
\draw (2.5,-1) node[anchor= south east] {$x$};
\draw (3.5,1) node[anchor=  west] {$u$} -- (3.5,-1.2);
\draw (3.5,-1) node[anchor= south west] {$v$};
\draw (2.5,3)[->]  .. controls (2.5,2) and (3.5,2) .. (3.5,1) ;
\draw (3.5,3)[->]  .. controls (3.5,2) and (2.5,2) .. (2.5,1) ;
\draw (3.5,0) node[morphism] {$f$};
\end{tikzpicture}
\qquad
\begin{tikzpicture}[line width=2pt]
\fill[color=blue!20] (-1,-1.2) rectangle (1,3.2);
\draw (-0.5,-1) node[anchor=south east] {$v$} -- (-0.5,-1.2);
\draw (0.5,-1) node[anchor=south west] {$x$} -- (0.5,-1.2);
\draw (-0.5,3.2)[->] -- (-0.5,1) node[anchor=  east] {$x$};
\draw (0.5,3.2)[->] -- (0.5,1) node[anchor=  west] {$v$};
\draw (-0.5,3.2)[->] -- (-0.5,3) node[anchor=north east] {$x$};
\draw (0.5,3.2)[->] -- (0.5,3)  node[anchor=north west] {$u$};
\draw (-0.5,1)[->]   .. controls (-0.5,0) and (0.5,0) .. (0.5,-1) ;
\draw (0.5,1)[->]  .. controls (0.5,0) and (-0.5,0) .. (-0.5,-1) ;
\draw (0.5,2) node[morphism] {$f$};
\draw (1.5,1) node {$=$};

\fill[color=blue!20] (2,-1.2) rectangle (4,3.2);
\draw (2.5,3.2)[->] -- (2.5,3) node[anchor=north east] {$x$};
\draw (3.5,3.2)[->] -- (3.5,3)  node[anchor=north west] {$u$};
\draw (2.5,1) -- (2.5,-1.2);
\draw (3.5,1) -- (3.5,-1.2);
\draw (2.5,-1) node[anchor = south east] {$v$};
\draw (3.5,-1) node[anchor = south west] {$x$};
\draw (2.5,3)[->]  .. controls (2.5,2) and (3.5,2) .. (3.5,1)
node[anchor=  west] {$x$};
\draw (3.5,3)[->]  .. controls (3.5,2) and (2.5,2) .. (2.5,1)
node[anchor=  east] {$u$};
\draw (2.5,0) node[morphism] {$f$};
\end{tikzpicture}
\end{center}
\caption{Relations \eqref{a:slide1} and
  \eqref{a:slide2}}\label{fig:rels-slide}
\end{figure}
\begin{figure}
\begin{center}
\begin{tikzpicture}[line width=2pt]
\fill[color=blue!20] (-1.75,-1.2) rectangle (1.75,1.2);
\draw (-0.5,-1) -- (-0.5,-1.2);
\draw (0.5,-1) -- (0.5,-1.2);
\draw (-0.5,1.2)[->] -- (-0.5,1);
\draw (0.5,1.2)[->] -- (0.5,1);
\draw (-0.5,1)[->] node[anchor=north east] {$x\otimes y$}
.. controls (-0.5,0) and (0.5,0) .. (0.5,-1) node[anchor=south west]
{$x\otimes y$};
\draw (0.5,1)[->] node[anchor=north west] {$z$}  .. controls (0.5,0)
and (-0.5,0) .. (-0.5,-1) node[anchor=south east] {$z$};
\draw (2.1,0) node {$=$};
\begin{scope}[xshift=3.5cm]
\fill[color=blue!20] (-1,-1.2) rectangle (2,1.2);
\draw (-0.5,-1) -- (-0.5,-1.2);
\draw (0.5,-1) -- (0.5,-1.2);
\draw (1.5,-1) -- (1.5,-1.2);
\draw (-0.5,1.2)[->] -- (-0.5,1);
\draw (0.5,1.2)[->] -- (0.5,1);
\draw (1.5,1.2)[->] -- (1.5,1);
\draw (-0.5,1)[->] node[anchor=north east] {$x$}  .. controls
(-0.5,0) and (0.5,0) .. (0.5,-1) node[anchor=south west] {$x$};
\draw (0.5,1)[->] node[anchor=north west] {$y$}  .. controls (0.5,0)
and (1.5,0) .. (1.5,-1) node[anchor=south west] {$y$};
\draw (1.5,1)[->] node[anchor=north west] {$z$}  .. controls (1.5,0)
and (-0.5,0) .. (-0.5,-1) node[anchor=south east] {$z$};
\end{scope}
\end{tikzpicture}
\qquad
\begin{tikzpicture}[line width=2pt]
\fill[color=blue!20] (-1.75,-1.2) rectangle (1.75,1.2);
\draw (-0.5,-1) -- (-0.5,-1.2);
\draw (0.5,-1) -- (0.5,-1.2);
\draw (-0.5,1.2)[->] -- (-0.5,1);
\draw (0.5,1.2)[->] -- (0.5,1);
\draw (-0.5,1)[->] node[anchor=north east] {$x$}  .. controls
(-0.5,0) and (0.5,0) .. (0.5,-1) node[anchor=south west] {$x$};
\draw (0.5,1)[->] node[anchor=north west] {$y\otimes z$}  .. controls
(0.5,0) and (-0.5,0) .. (-0.5,-1) node[anchor=south east] {$y\otimes
z$};
\draw (2.1,0) node {$=$};
\begin{scope}[xshift=3.5cm]
\fill[color=blue!20] (-1,-1.2) rectangle (2,1.2);
\draw (-0.5,-1) -- (-0.5,-1.2);
\draw (0.5,-1) -- (0.5,-1.2);
\draw (1.5,-1) -- (1.5,-1.2);
\draw (-0.5,1.2)[->] -- (-0.5,1);
\draw (0.5,1.2)[->] -- (0.5,1);
\draw (1.5,1.2)[->] -- (1.5,1);
\draw (-0.5,1)[->] node[anchor=north east] {$x$}  .. controls
(-0.5,0) and (1.5,0) .. (1.5,-1) node[anchor=south west] {$x$};
\draw (0.5,1)[->] node[anchor=north east] {$y$}  .. controls (0.5,0)
and (-0.5,0) .. (-0.5,-1) node[anchor=south west] {$y$};
\draw (1.5,1)[->] node[anchor=north west] {$z$}  .. controls (1.5,0)
and (0.5,0) .. (0.5,-1) node[anchor=south east] {$z$};
\end{scope}
\end{tikzpicture}
\end{center}
\caption{Relations \eqref{a:tensor1} and \eqref{a:tensor2}}\label{fig:rels-tensor}
\end{figure}
Note that applying the relations \eqref{a:slide1}, \eqref{a:slide2}
to $f=\alpha$
give the braid relations for the morphisms $\alpha$.

Braided monoidal categories are only mentioned briefly in this paper.
The definition is given by weakening the definition of a symmetric
monoidal category; replacing the symmetric groups by the braid groups
and replacing the string diagrams by braided versions.

\subsection{Duality}
Duals in monoidal categories are based on dual vector spaces. A vector space
has a dual if and only if it is finite dimensional. In this paper all
representations are finite dimensional and have a dual.
\begin{defn}
Let $x$ and $y$ be objects in a monoidal category. Then $x$ is a
\emph{left dual}
of $y$ and $y$ is a \emph{right dual} of $x$ means that we are given
evaluation and
coevaluation morphisms
$I\rightarrow x\otimes y$ and $y\otimes x\rightarrow I$ such that
both of the
following composites are identity morphisms
\begin{align*}
 x\rightarrow I\otimes x\rightarrow x\otimes y\otimes x\rightarrow
x\otimes I\rightarrow x \label{s:trl} \\
 y\rightarrow y\otimes I\rightarrow y\otimes x\otimes y\rightarrow
I\otimes y\rightarrow y. 
\end{align*}
The object $x$ is \emph{dual} to $y$ if it is both a left dual and a
right dual.
\end{defn}

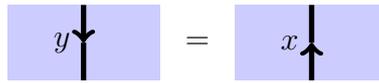
\begin{figure}
\begin{center}
 \begin{tikzpicture}[line width=2pt]
\fill[color=blue!20] (-1,-0.5) rectangle (1,0.5);
\draw (0,0.5)[->] -- (0,0) node[anchor=east] {$y$};
\draw (0,0) -- (0,-0.5);
\draw (1.5,0) node {$=$};
\fill[color=blue!20] (2,-0.5) rectangle (4,0.5);
\draw (3,-0.5)[->] -- (3,0) node[anchor=east] {$x$};
\draw (3,0) -- (3,0.5);
 \end{tikzpicture}
 \end{center}
\caption{Convention for $x$ being a dual of $y$}\label{fig:dual}
\end{figure}

In the string diagrams we adopt the convention depicted in
Figure~\ref{fig:dual}.
Using this convention and omitting the edge labelled $I$ the string diagrams
for the evaluation and coevaluation morphisms are:
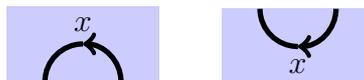
\begin{equation*}
\begin{tikzpicture}[line width=2pt]
\fill[color=blue!20] (-1,-1) rectangle (1,0);
\draw (0.5,-1)[->] arc (0:90:0.5);
\draw (0,-0.5) node[anchor=south] {$x$} arc (90:180:0.5);
\end{tikzpicture}
\qquad
\begin{tikzpicture}[line width=2pt]
\fill[color=blue!20] (-1,-1) rectangle (1,0);
\draw (0.5,0)[->] arc (0:-90:0.5);
\draw (0,-0.5) node[anchor=north] {$x$} arc (-90:-180:0.5);
\end{tikzpicture}
\end{equation*}

\begin{figure}
\begin{center}
\begin{tikzpicture}[line width=2pt]
\fill[color=blue!20] (-1.25,-1) rectangle (1.25,1);
\draw[->] (1,1) -- (1,0) arc(0:-90:0.5);
\draw[->] (0.5,-0.5) node[anchor=north] {$x$} arc(-90:-180:0.5)
arc(0:90:0.5);
\draw (-0.5,0.5) node[anchor=south] {$x$} arc(90:180:0.5) -- (-1,-1);
\draw (1.5,0) node {$=$};
\fill[color=blue!20] (1.75,-1) rectangle (2.5,1);
\draw[->] (2,1) -- (2,0);
\draw (2,0) node[anchor=west] {$x$} -- (2,-1);
\end{tikzpicture}
\qquad\qquad
\begin{tikzpicture}[line width=2pt]
\fill[color=blue!20] (-1.25,-1) rectangle (1.25,1);
\draw[->] (1,-1) -- (1,0) arc(0:90:0.5);
\draw[->] (0.5,0.5) node[anchor=south] {$x$} arc(90:180:0.5) arc(0:-90:0.5);
\draw (-0.5,-0.5) node[anchor=north] {$x$} arc(-90:-180:0.5) -- (-1,1);
\draw (1.5,0) node {$=$};
\fill[color=blue!20] (1.75,-1) rectangle (2.5,1);
\draw[->] (2,-1) -- (2,0);
\draw (2,0) node[anchor=west] {$x$} -- (2,1);
\end{tikzpicture}
\end{center}
\caption{Relations for duals}\label{fig:cupcap}
\end{figure}
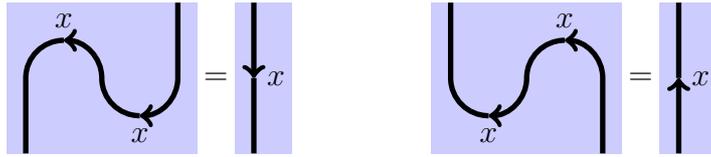

The string diagrams for the conditions on these morphisms are shown
in Figure~\ref{fig:cupcap}.

It follows that there are natural homomorphisms
\begin{equation*}
  \Hom(y\otimes u,v)\rightarrow \Hom(u,x\otimes v)
\end{equation*}
sending $\phi\colon y\otimes u\rightarrow v$ to the composite
\begin{equation*}
  u\rightarrow I\otimes u\rightarrow x\otimes y\otimes u\rightarrow
  x\otimes v.
\end{equation*}
Similarly there are natural homomorphisms
\begin{equation*}
  \Hom(u,x\otimes v)\rightarrow \Hom(y\otimes u,v).
\end{equation*}
sending $\phi\colon u\rightarrow x\otimes v$ to the composite
\begin{equation*}
  y\otimes u\rightarrow y\otimes x\otimes v\rightarrow I\otimes
  v\rightarrow  v.
\end{equation*}
It follows from the relations in Figure~\ref{fig:cupcap} that these
are inverse isomorphisms.

\begin{ex} Take a category $\sC$ and consider the category whose objects are
functors $F\colon \sC \rightarrow \sC$ and whose morphisms are natural
transformations. This is a monoidal category with $\otimes$ given by
composition.
Then $F$ is left dual to $G$ in this monoidal category is equivalent to
$F$ is left adjoint to $G$.
\end{ex}

\begin{ex} Let $K$ be a field. The category of vector spaces over $K$ is
a symmetric monoidal category. A vector space has a dual if and only if
it is finite dimensional. More generally, Let $K$ be a commutative ring.
The category of $K$-modules is a symmetric monoidal category.
A $K$-module has a dual if and only if it is finitely generated and
projective.
\end{ex}

\subsection{Rotation}
Let $x$ be an object in a monoidal category. Then $\Hom(I,\otimes^r
x)$ is the set of \emph{invariant tensors}.  In this section we
discuss two constructions of an action of the cyclic group on this
set each of which depends on additional structure. One construction,
which we call rotation, requires that $x$ has a dual.  The other
construction assumes that the category is symmetric monoidal. The
main result in this section compares these two actions in the
situation when both are defined.

Let $x$ be an object of a symmetric monoidal category.  Then
$\otimes^r x$ has a natural action of $\fS_r$ where the action of the
Coxeter generator $s_i$ is given by the isomorphism
$(\otimes^{i-1}\id_x)\otimes\alpha(x,x)\otimes
(\otimes^{r-i-1}\id_x)$.  This induces an action of $\fS_r$ on the
set $\Hom(y,\otimes^r x)$ for any $y$. In particular for $y=I$ this
is an action of $\fS_r$ on the invariant tensors.

Assume $x$ has a dual $x^\ast$. Then, for any $y$, there is a natural
map $\Hom(I,y\otimes x)\rightarrow \Hom(I,x\otimes y)$ sending
$f\colon I\rightarrow y\otimes x$ to the composite
\begin{equation*}
 I\rightarrow x\otimes x^*\rightarrow x\otimes I\otimes x^*\rightarrow
x\otimes y\otimes x\otimes x^\ast\rightarrow x\otimes y\otimes I
\rightarrow x\otimes y
\end{equation*}

Taking $y=\otimes^{r-1}x$ gives the \emph{rotation map}
$\otimes^rx\rightarrow \otimes^rx$. Then, for any $x$, the $r$-th
power of the rotation map is the identity.

\begin{prop}\label{lem:rot} Let $x$ be an object of symmetric monoidal
category with dual $x^\ast$ such that the condition in Figure \ref{fig:con}
is satisfied.
\begin{figure}
\setlength{\tabcolsep}{0.6cm}
\begin{center}
\begin{tabular}{cc}
\begin{tikzpicture}[line width=2pt]
\fill[color=blue!20] (-1,-1) rectangle (1,1);
\draw[->] (0.5,-1) .. controls (0.5,-0.5) and (-0.5,-0.5) .. (-0.5,0) arc(180:90:0.5);
\draw (0,0.5) node[anchor=south] {$x$} arc(90:0:0.5) .. controls (0.5,-0.5) and (-0.5,-0.5) .. (-0.5,-1);
\draw (1.5,0) node {$=$};
\fill[color=blue!20] (2,-1) rectangle (4,1);
\draw[->] (2.5,-1) -- (2.5,0) arc(180:90:0.5);
\draw (3,0.5) node[anchor=south] {$x$} arc(90:0:0.5) -- (3.5,-1);
\end{tikzpicture} &
\begin{tikzpicture}[line width=2pt]
\fill[color=blue!20] (-1,-1) rectangle (1,1);
\draw[->] (-0.5,-1) .. controls (-0.5,-0.5) and (0.5,-0.5) .. (0.5,0)
arc(0:90:0.5);
\draw (0,0.5) node[anchor=south] {$x$} arc(90:180:0.5) .. controls
(-0.5,-0.5) and (0.5,-0.5) .. (0.5,-1);
\draw (1.5,0) node {$=$};
\fill[color=blue!20] (2,-1) rectangle (4,1);
\draw[->] (3.5,-1) -- (3.5,0) arc(0:90:0.5);
\draw (3,0.5) node[anchor=south] {$x$} arc(90:180:0.5) -- (2.5,-1);
\end{tikzpicture} \\*[1cm]
\begin{tikzpicture}[line width=2pt]
\fill[color=blue!20] (-1,-1) rectangle (1,1);
\draw[->] (-0.5,1) .. controls (-0.5,0.5) and (0.5,0.5) .. (0.5,0)
arc(0:-90:0.5);
\draw (0,-0.5) node[anchor=north] {$x$} arc(-90:-180:0.5) .. controls
(-0.5,0.5) and (0.5,0.5) .. (0.5,1);
\draw (1.5,0) node {$=$};
\fill[color=blue!20] (2,-1) rectangle (4,1);
\draw[->] (2.5,1) -- (2.5,0) arc(180:270:0.5);
\draw (3,-0.5) node[anchor=north] {$x$} arc(-90:0:0.5) -- (3.5,1);
\end{tikzpicture} &
\begin{tikzpicture}[line width=2pt]
\fill[color=blue!20] (-1,-1) rectangle (1,1);
\draw[->] (0.5,1) .. controls (0.5,0.5) and (-0.5,0.5) .. (-0.5,0)
arc(180:270:0.5);
\draw (0,-0.5) node[anchor=north] {$x$} arc(-90:0:0.5) .. controls
(0.5,0.5) and (-0.5,0.5) .. (-0.5,1);
\draw (1.5,0) node {$=$};
\fill[color=blue!20] (2,-1) rectangle (4,1);
\draw[->] (2.5,1) -- (2.5,0) arc(180:270:0.5);
\draw (3,-0.5) node[anchor=north] {$x$} arc(-90:0:0.5) -- (3.5,1);
\end{tikzpicture}
\end{tabular}
\end{center}
\caption{Conditions on a representation}\label{fig:con}
\end{figure}
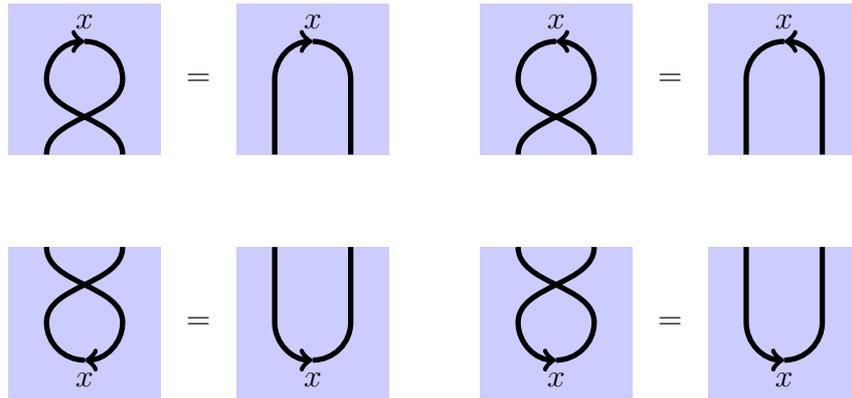
Then the action of the long cycle is given by rotation.
\end{prop}
It is clear that all the conditions on the string diagrams can be
interpreted as sliding the strings in the diagrams.  The converse is
that if the strings in one diagram can be slid to obtain another
diagram then the two morphisms are equal.

\begin{proof}
  In the following diagrams the half circle represents an invariant
  tensor.  The first diagram shows the map given by the symmetry and
  the final diagram shows the rotation map.
  
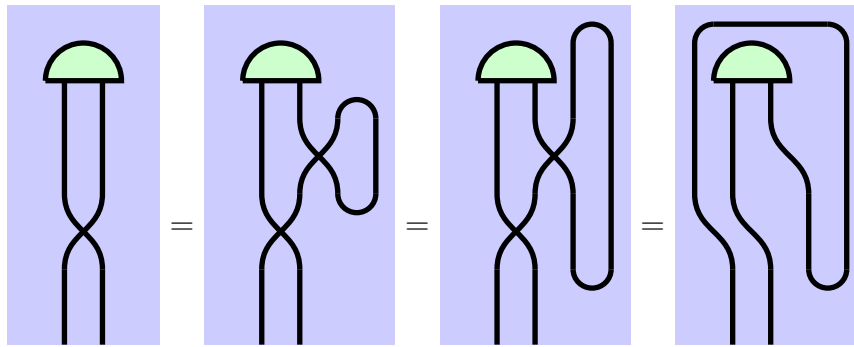
\begin{center}
    \begin{tikzpicture}[line width = 2pt]
      \fill[color=blue!20] (0,-1) rectangle  (2,3.5);
      \fill[morphism] (0.5,2.5) -- (1.5,2.5) arc (0:180:0.5);
      \draw (0.75,2.5) -- (0.75,1) .. controls (0.75,0.5) and (1.25,0.5) .. (1.25,0);
      \draw (1.25,2.5) -- (1.25,1) .. controls (1.25,0.5) and (0.75,0.5) .. (0.75,0);
      \draw (0.75,0) -- (0.75,-1);
      \draw (1.25,0) -- (1.25,-1);
    \end{tikzpicture}
    \raisebox{1.5cm}{$=$}
    \begin{tikzpicture}[line width = 2pt]
      \fill[color=blue!20] (0,-1) rectangle  (2.5,3.5);
      \fill[morphism] (0.5,2.5) -- (1.5,2.5) arc (0:180:0.5);
      \draw (0.75,2.5) -- (0.75,1) .. controls (0.75,0.5) and  (1.25,0.5) .. (1.25,0);
      \draw (1.25,2.5) -- (1.25,2) .. controls (1.25,1.5) and (1.75,1.5) .. (1.75,1);
      \draw (1.75,2) .. controls (1.75,1.5) and (1.25,1.5) .. (1.25,1);
      \draw (1.25,1) .. controls (1.25,0.5) and (0.75,0.5) .. (0.75,0);
      \draw (1.75,2) arc(180:0:0.25);
      \draw (2.25,2) -- (2.25,1);
      \draw (1.75,1) arc(-180:0:0.25);
      \draw (0.75,0) -- (0.75,-1);
      \draw (1.25,0) -- (1.25,-1);
    \end{tikzpicture}
    \raisebox{1.5cm}{$=$}
    \begin{tikzpicture}[line width = 2pt]
      \fill[color=blue!20] (0,-1) rectangle  (2.5,3.5);
      \fill[morphism] (0.5,2.5) -- (1.5,2.5) arc (0:180:0.5);
      \draw (0.75,2.5) -- (0.75,1) .. controls (0.75,0.5) and  (1.25,0.5) .. (1.25,0);
      \draw (1.25,2.5) -- (1.25,2) .. controls (1.25,1.5) and (1.75,1.5) .. (1.75,1);
      \draw (1.75,2) .. controls (1.75,1.5) and (1.25,1.5) .. (1.25,1);
      \draw (1.25,1) .. controls (1.25,0.5) and (0.75,0.5) .. (0.75,0);
      \draw (1.75,3) arc(180:0:0.25);
      \draw (2.25,3) -- (2.25,0);
      \draw (1.75,3) -- (1.75,2);
      \draw (1.75,1) -- (1.75,0);
      \draw (1.75,0) arc(-180:0:0.25);
      \draw (0.75,0) -- (0.75,-1);
      \draw (1.25,0) -- (1.25,-1);
    \end{tikzpicture}
    \raisebox{1.5cm}{$=$}
    \begin{tikzpicture}[line width = 2pt]
      \fill[color=blue!20] (0,-1) rectangle  (2.5,3.5);
      \fill[morphism] (0.5,2.5) -- (1.5,2.5) arc (0:180:0.5);
      \draw (0.75,2.5) -- (0.75,1) .. controls (0.75,0.5) and  (1.25,0.5) .. (1.25,0);
      \draw (1.25,2.5) -- (1.25,2) .. controls (1.25,1.5) and (1.75,1.5) .. (1.75,1);
      \draw (0.25,3) -- (0.25,1) .. controls (0.25,0.5) and (0.75,0.5) .. (0.75,0);
      \draw (2.25,3) arc(0:90:0.25);
      \draw (2,3.25) -- (0.5,3.25);
      \draw (0.5,3.25) arc(90:180:0.25);
      \draw (2.25,3) -- (2.25,0);
      \draw (1.75,1) -- (1.75,0);
      \draw (1.75,0) arc(-180:0:0.25);
      \draw (0.75,0) -- (0.75,-1);
      \draw (1.25,0) -- (1.25,-1);
    \end{tikzpicture}
\end{center}
The first equation is a consequence of the relations for duals in
Figure~\ref{fig:cupcap}, relation~\eqref{a:slide1} in
Figure~\ref{fig:rels-slide} (with $f$ being the coevaluation) and the
condition in Figure~\ref{fig:con}.  The second equation is the
condition that two tensors commute.  The third equation is
an application of~\eqref{a:slide1} and~\eqref{a:tensor1}.

For illustration, we demonstrate the first equation in more detail:
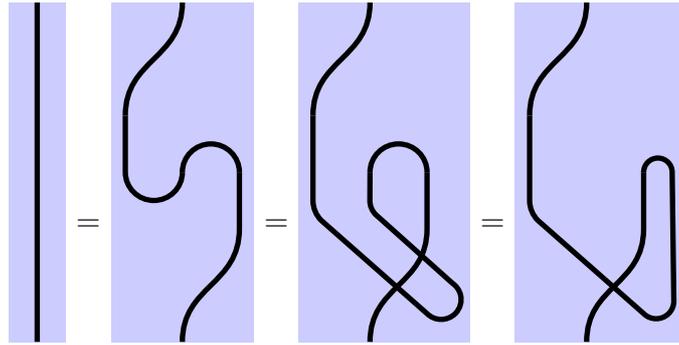
\begin{center}
  \begin{tikzpicture}[line width=2pt, scale=0.75]
    \fill[color=blue!20] (-0.5,-3) rectangle (0.5,3);
    \draw (0,-3) -- (0,3);
  \end{tikzpicture}
  \raisebox{1.5cm}{$=$}
  \begin{tikzpicture}[line width=2pt, scale=0.75]
    \fill[color=blue!20] (-1.25,-3) rectangle (1.25,3);
    \draw (0,3) .. controls (0,2) and (-1,2) .. (-1,1);
    \draw (-1, 1) -- (-1, 0) arc(180:360:0.5);
    \draw (1,0) arc(0:180:0.5);
    \draw (1, 0) -- (1, -1) .. controls (1,-2) and (0,-2) .. (0,-3);
  \end{tikzpicture}
  \raisebox{1.5cm}{$=$}
  \begin{tikzpicture}[line width=2pt, scale=0.75]
    \fill[color=blue!20] (-1.25,-3) rectangle (1.75,3);
    \draw (0,3) .. controls (0,2) and (-1,2) .. (-1,1);
    \draw (-1, 1) -- (-1, -0.5) arc(180:225:0.5) -- (1, -2.5)
    arc(225:415:0.35) -- (0.1, -0.75) arc(225:180:0.35) -- (0, 0);
    \draw (1,0) arc(0:180:0.5);
    \draw (1, 0) -- (1, -1) .. controls (1,-2) and (0,-2) .. (0,-3);
  \end{tikzpicture}
  \raisebox{1.5cm}{$=$}
  \begin{tikzpicture}[line width=2pt, scale=0.75]
    \fill[color=blue!20] (-1.25,-3) rectangle (1.75,3);
    \draw (0,3) .. controls (0,2) and (-1,2) .. (-1,1);
    \draw (-1, 1) -- (-1, -0.5) arc(180:225:0.5) -- (1, -2.5)
    arc(225:360:0.31) -- (1.5, 0) arc(0:180:0.25);
    \draw (1, 0) -- (1, -1) .. controls (1,-2) and (0,-2) .. (0,-3);
  \end{tikzpicture}
\end{center}
\end{proof}

The strategy then for constructing examples of the cyclic sieving phenomenon
is the following. Let $x$ be an object in a symmetric monoidal category
with a dual $x^\ast$. The combinatorial structure is a basis of $\Hom(I,\otimes^n x)$
which is preserved by rotation. Then the aim is to apply Theorem \ref{thm:fund}
to obtain an example of the cyclic sieving phenomenon.

\subsection{Coboundary categories}
The representations of a quantum group are famously not a symmetric monoidal
category. Instead it has two weaker structures both introduced by Drinfel$'$d.
One is a braiding which has applications to knot theory. The other is a
coboundary structure which has received less attention.

\begin{defn} A coboundary category is a monoidal category with natural isomorphism
$\sigma_{A,B}\colon A\otimes B \rightarrow B\otimes A$ such that 
\[ \sigma_{A,B}\circ \sigma_{B,A}=1 \]
for all $A$, $B$ and such that for all $A,B,C$ the following diagram commutes
\begin{equation} \begin{CD} A\otimes B\otimes C @>{1\otimes\sigma_{B,C}}>>  A\otimes C\otimes B \\
@V{\sigma_{A,B}\otimes 1}VV @VV{\sigma_{A,C\otimes B}}V \\
 B\otimes A\otimes C @>>{\sigma_{B\otimes A,C}}>  C\otimes B\otimes A
\end{CD} \end{equation} 
\end{defn}

The \emph{cactus groups} are a sequence of groups, $\{\cactus_r:r\ge 0\}$, which share many of the
properties of the braid groups. This analogy is developed in \cite{MR1718078} and \cite{MR2448884}.
There are surjective homomorphisms $\cactus_r\rightarrow \fS_r$
and the kernels are the pure cactus groups.
\begin{defn} The group $\cactus_r$ is generated by elements $\langle u:v:w\rangle$
for $0\le u<v<w<r$. 
Defining relations are
\begin{align}
\langle u:v:w\rangle \langle u':v':w'\rangle &= \langle u':v':w'\rangle \langle u:v:w\rangle \\
&\qquad \text{if $u'>w$ or $w'<u$} \nonumber \\
\langle u:v:w\rangle &\langle u:u+w-v:w\rangle =1 \\
\langle x:y:z\rangle\langle u:v:w\rangle &= \langle u:v:w\rangle \langle x+w-v:y+w-v:z+w-v\rangle \\
&\qquad \text{if $u<x$ and $z\le v$} \nonumber \\
\langle x:y:z\rangle\langle u:v:w\rangle &= \langle u:v:w\rangle \langle x+u-v:y+u-v:z+u-v\rangle \\
&\qquad \text{if $v<x$ and $z\le w$} \nonumber \\
\langle v:w:x\rangle\langle u:v:x\rangle &= \langle u:v:w\rangle \langle u:w:x\rangle
\end{align}
\end{defn}
These are the commutation relations, the symmetry relations, the two naturality relations
and the coboundary relations.

The homomorphism to the symmetric group is given by
\[ \langle u:v:w\rangle(k) = \begin{cases}
k+w-v &\text{if $u<k\le v$} \\
k+u-v &\text{if $v<k\le w$} \\
k  & \text{otherwise}
\end{cases} \]

\begin{lemma} Let $x$ be an object in a coboundary category. Define a natural action
of the generator $\langle u:v:w\rangle$ on $\otimes^rx$ by
\begin{equation*}
\langle u:v:w\rangle = (\otimes^{u-1}\id_x) \otimes \sigma_{\otimes^{v-u}x,\otimes^{w-v}x} \otimes (\otimes^{r-w}\id_x)
\end{equation*}
Then these satisfy the defining relations and so give 
a natural action of $\cactus_r$ on $\otimes^rx$.
\end{lemma}

\section{Invariant tensors}
The aim of this article is to construct examples of the cyclic sieving phenomenon
by applying Theorem \ref{thm:fund} to spaces of invariant tensors. More specifically,
let $\fg$ be a semisimple complex Lie algebra with enveloping algebra $U(\fg)$. Denote the set of positive weights by $P_+$; this parametrises the set of isomorphism classes of irreducible representations. Let $M$ be a finite dimensional
representation. Then $\otimes^rM$ is also a representation and we denote by
$N(\otimes^rM)$ the subspace of invariant tensors. This space has a natural action
of $\fS_r$. In order to apply Theorem \ref{thm:fund} we need to determine the
Frobenius character and to construct a basis invariant under the long cycle.

In this article we assume, for simplicity, that the representation $M$ is irreducible.
The following Lemma allows this assumption to be dropped.
\begin{lemma} For all $M$ and $M'$, 
\begin{equation}
\bch N(\otimes^r(M\oplus M')) = 
\sum_{s=0}^r \bch N(\otimes^sM)\,\bch N(\otimes^{r-s}M')
\end{equation}	
\end{lemma}

In Proposition \ref{lem:rot} there is a restriction on the representation $M$.
For $M$ irreducible there are three cases. One case is $M\not\cong M^*$.
The alternative $M\cong M^*$ gives two cases; $M$ has a nondegenerate symmetric
inner product and $M$ has a symplectic form.

If $M$ has a nondegenerate symmetric
inner product then $M$ satisfies the condition in Figure \ref{fig:con}. If $M\not\cong M^*$ then there are two evaluation maps and two coevaluation maps and these maps can always
be chosen to satisfy the condition in Figure \ref{fig:con}.

If $M$ has a symplectic form then we use an alternative sign convention. A vector space with
a symplectic form can also be regarded as an odd super vector space
with a symmetric inner product. So we regard $M$ as an odd super
representation; this satisfies the condition in Figure \ref{fig:con}.
Then $N(\otimes^rM)$ as a representation of $\fS_r$ is tensored with
the sign representation and the involution $\omega$ is applied to the
Frobenius character.

Next we show that the symmetric functions  $\bch N(\otimes^rM)$ are determined by
the characters of the Schur functors (aka plethysms) applied to $M$.
Since the operations on characters associated to Schur functors have been implemented 
in the LiE package, \cite{vLCL92} this means that the Frobenius characters can be computed
in small examples. This package is no longer supported
but it has been incorporated into Magma, \cite{MR1484478} and Sage, \cite{sage}.
The calculations in this paper used these two packages.

Each representation, $V$, of $\fS_r$ gives a polynomial functor on representations of $\fg$
by $M\mapsto (\otimes^rM)\otimes_{\bC\fS_r} V$. For $\lambda\vdash r$ let $S(\lambda)$
be the corresponding irreducible representation of $\fS_r$ and $\bS^\lambda$ the
associated polynomial functor. The functor $\bS^\lambda$ is known as a Schur functor.

The decomposition of $\otimes^r M$ as a representation of $\fS(r)\times \GL(M)$ is
\[ \otimes^r M \cong \bigoplus_{\lambda \vdash r} S(\lambda)\otimes \bS^\lambda(M) \]
Now regard $\bS^\lambda(M)$ as a $\fg$-module and take the decomposition into irreducible
representations
\[ \bS^\lambda(M)\cong \bigoplus_{\varpi\in P_+} A(\lambda,\varpi)\otimes V(\varpi) \]
where $A(\lambda,\varpi)$ is a vector space.

Then the decomposition of
$\otimes^r M$ as a representation of $\bC \fS(r)\otimes U(\fg)$ is
\[  \otimes^r M \cong \bigoplus_{\substack{\lambda\vdash r\\ \varpi\in P_+}} 
 A(\lambda,\varpi)\otimes S(\lambda)\otimes V(\varpi) \]

Define $a(\lambda,\varpi)=\dim A(\lambda,\varpi)$. Then the Frobenius character
of the isotypical subspace of $\otimes^rM$ associated to $\varpi\in P_+$ is
$\sum_{\lambda\vdash r}a(\lambda,\varpi)s_\lambda$. In particular, taking $\varpi=0$,
\begin{equation*}
 \bch N(\otimes^rM) = \sum_{\lambda\vdash r} a(\lambda,0)s_\lambda
\end{equation*}

\section{Quantum groups}\label{sec:qg}
In this section we give a summary of the basic results on the representation
theory of the Drinfel$'$d-Jimbo quantised enveloping algebra of a semisimple
Lie algebra. Our standard references are \cite[Part 1]{MR1227098} and
\cite[Chapters 4 \& 5]{MR1359532}.

Fix a finite type Cartan matrix, $\sC$. The structures associated to $\sC$
that we will discuss in this section are the semisimple complex Lie algebra $\fg(\sC)$,
its universal enveloping algebra $U(\sC)$ (considered as a Hopf algebra over $\bQ$);
the Jimbo quantised enveloping algebra $U_q(\sC)$ (considered as a Hopf algebra over $\bQ(q)$);
and the categories of finite dimensional representations of these Hopf algebras.

\subsection{Presentation}
The quantized enveloping algebra $U_q=U_q(\mathfrak{g})$ is the associative
algebra (with one) over $\bQ(q)$ generated by $F_{\beta}$, $E_{\beta}$ for $\beta\in\Pi$
and $K_\alpha$ for $\alpha\in Q$. The following are the defining
relations which hold for
all $\alpha,\alpha'\in Q$ and all $\beta,\beta'\in\Pi$:
\begin{align*}
K_0&=1 \\
K_\alpha K_{\alpha'} &= K_{\alpha+\alpha'} \\
K_{\alpha} E_{\beta} &=  q^{(\alpha,\beta)}E_{\beta}K_{\alpha}\\
F_{\beta} K_{\alpha} &= q^{-(\alpha,\beta)}K_{\alpha} F_{\beta}\\
E_{\beta} F_{\beta'} - F_{\beta'}E_{\beta}  &= \delta_{\beta,\beta'}
\frac{K_{\beta}-K_{\beta}^{-1}}{q_{\beta}-q_{\beta}^{-1}}
\end{align*}
together with the quantum Serre relations which are omitted.

This is a Hopf algebra with coproduct, counit, and antipode. The counit gives
the trivial representation, the coproduct gives the tensor product and the
antipode gives the dual representation.

The counit is
\begin{equation*}
K_{\alpha}\mapsto 1\qquad E_{\alpha}\mapsto 0\qquad F_{\alpha}\mapsto 0
\end{equation*}

There are many possible ways to define a comultiplication
$\Delta : U_q\to U_q\otimes U_q$. One choice is:
\begin{align*}
\Delta(E_{\alpha}) &= E_{\alpha}\otimes 1 + K_{\alpha}\otimes E_{\alpha}\\
\Delta(F_{\alpha}) &= F_{\alpha}\otimes K_{\alpha}^{-1} +1\otimes F_{\alpha}\\
\Delta(K_{\alpha}) &= K_{\alpha}\otimes K_{\alpha}
\end{align*}
The associated antipode is
\begin{equation*}
S(K_\alpha)=K_\alpha^{-1}\qquad S(E_\alpha)=-K_\alpha^{-1}E_\alpha
\qquad S(F_\alpha)=-F_\alpha K_\alpha
\end{equation*}

\begin{ex} The first example, both historically and pedagogically,
	is the case $\SL(2)$. This has rank one and the Cartan matrix is $(2)$.
	The algebra $U_q(\SL(2))$ is generated by elements $K^\pm 1$, $E$, $F$
	and the defining relations are
	\begin{align*}
	KK^{-1} &= 1 & EK &= qKE \\
	K^{-1}K &= 1 & FK &= q^{-1}KF \\
	\end{align*}
	\begin{equation*}
	EF-FE = \frac{K-K^{-1}}{q-q^{-1}}
	\end{equation*}
\end{ex}

The integral form is an analogue of the Chevalley integral form and was introduced in \cite[\S 4]{MR954661}.
Recall that the $q$-factorial is defined by
\begin{equation*}
[k]=\frac{q^k-q^{-k}}{q-q^{-1}} \quad [k]!=[k]\,[k-1]\, \dotsc \,[1]
\end{equation*}
\begin{defn}
	The integral form is the $\bZ[q,q^{-1}]$ subalgebra of $U_q(\sC)$ generated by the generators $K_\alpha$ and the quantum divided powers
	\begin{equation}
	E_\alpha^{(k)}=\frac{E_\alpha^k}{[k]!}\qquad
	F_\alpha^{(k)}=\frac{F_\alpha^k}{[k]!}
	\end{equation}
\end{defn}

The Lusztig involution was introduced in \cite{MR1182165}. Fix a finite type Cartan matrix $\sC$. Let $\Pi$ be the set of simple roots, $w_0$ be the longest element in the Weyl group and let $P_+$ be the dominant weights.

Let $\theta\colon \Pi\rightarrow \Pi$ be the Dynkin diagram automorphism such that $\theta(\alpha)=-w_0\alpha$.

\begin{defn} The Lusztig involution is an $\bQ(q)$-algebra involution $\xi\colon U_q(\sC)\rightarrow U_q(\sC)$. It is defined on the generators by
	\begin{equation}
	\xi(E_\alpha)=F_{\theta(\alpha)}\quad
	\xi(F_\alpha)=E_{\theta(\alpha)}\quad
	\xi(K_\alpha)=K_\alpha^{-1}
	\end{equation}
\end{defn}

\subsection{Representations}
Let $U(\sC)-mod$ be the category of finite dimensional representations of
$U(\sC)$. The category $U(\sC)-mod$ is semisimple abelian and the simple objects
are the highest weight representations $L(\varpi)$ for $\varpi\in P_+$, the dominant
weights. It is often said that $U_q(\sC)$ is a deformation of $U(\sC)$ as a Hopf algebra.
This is inaccurate as $U_q(C)$ has more irreducible representations. However there is
a full subcategory, $U_q(\sC)-mod$, of the category of finite dimensional
representations of $U_q(\sC)$ such that the simple objects
are the highest weight representations $L_q(\varpi)$ for $\varpi\in P_+$
and which is a deformation of $U(\sC)-mod$.

\begin{ex}\label{ex:sl2} For each $n\ge 0$, there is a representation, $L(n)$, of $U_q(\SL(2))$ of dimension
$n+1$. Take the basis to be $v_0,v_1,\dotsc ,v_n$. Then the action of the generators
is given by $Kv_i=q^{n-2i}v_i$ and
\begin{equation}
Kv_i = q^{n-2i}v_i \qquad
Fv_i = [i+1]v_{i+1} \qquad
Ev_i = [n-i+1]v_{i-1}
\end{equation}
with the understanding that $v_{-1}=0$ and $v_{n+1}=0$.
\end{ex}

The categories $U(\sC)-mod$ and $U_q(\sC)-mod$ are both pivotal and so we also have
an action of the cyclic group $C_r$ on $\otimes^r V$ for any $V\in U(\sC)-mod$ and $r\ge 0$
where the action of the generator is given by rotation.

\begin{defn} The bar involution is the involution of $U_q(\sC)$ as a $\bQ$-algebra which is defined on the generators by
\begin{equation}
\overline{E_\alpha}=E_\alpha \quad
\overline{F_\alpha}=F_\alpha \quad
\overline{K_\alpha}=K_\alpha^{-1} \quad
\overline{q}=q^{-1}
\end{equation}
\end{defn}

Let $M\in U_q(\sC)-mod$. Then a Lusztig involution on $M$ is a $\bQ(q)$-linear map
$\xi_M\colon M\rightarrow M$ such that
\begin{equation}
\xi_M(um)=\xi(u)\,\xi_M(m)
\end{equation}
for $u\in U_q(\sC)$ and $m\in M$.

In Example \ref{ex:sl2}, which defines $L(n)$, the Lusztig involution is given by $\xi_{L(n)}(v_i)=v_{n-i}$.

The Kashiwara operators appear in the definition of a based module.
The definition is elementary but technical and is omitted.
In Example \ref{ex:sl2} the Kashiwara operators are given by
$\widetilde{F}v_i = v_{i+1}$ and $\widetilde{E}v_i = v_{i-1}$.

\subsection{Coboundary category}\label{sec:un}
There are two constructions which give $U_q(\sC)-mod$ the structure of
a coboundary monoidal category. The original construction by Drinfel$'$d
in \cite{MR802128} and \cite{MR934283} uses the unitarised $R$-matrix. Let $R$ be the $R$-matrix
then the unitarised $R$-matrix is $R(R^{\mathrm{op}}R)^{-1/2}$. This is defined using holomorphic functional calculus.


Let $\varpi\in P_+$ and denote the corresponding highest weight representation by $L_q(\varpi)$.
Then for $r>0$ there is an action of the cactus group $\cactus_r$ on $\otimes^r L_q(\varpi)\in U_q(\sC)-mod$. This gives a $\bQ(q)$-linear action on each isotypic subspace. Assume we are given the action of the Artin generators of the braid group $B_r$ on one of these vector spaces, then we can apply unitarisation to compute the action of each generator $\langle u:v:w\rangle\in\cactus_r$.

There is a standard way to lift a permutation to a positive braid. Take a reduced word for the permutation in the Coxeter generators and read this as a word in the Artin generators.

Take $\langle u:v:w\rangle\in\cactus_r$, map to a permutation and lift to a positive braid. This gives the braid
\begin{equation*}
(\sigma_{v}\dotsb \sigma_{u+1})(\sigma_{v+1}\dotsb \sigma_{u+2})\dotsb (\sigma_{w-1}\dotsb \sigma_{u-v+w})
\end{equation*}
Let $R$ be the matrix which represents this word and 
$R^{\mathrm{op}}$ be the matrix which represents the reversed word.

The matrix $R^{\mathrm{op}}R$ is semisimple and all eigenvalues are of the form $q^{2k}$ for $k\in\bZ$. This means it has a spectral decomposition as
\begin{equation*}
R^{\mathrm{op}}R = \sum_i q^{2k_i}E_i
\end{equation*}
where the $E_i$ are orthogonal idempotents.
Then $(R^{\mathrm{op}}R)^{-1/2}$ is the matrix
\begin{equation*}
(R^{\mathrm{op}}R)^{-1/2} = \sum_i q^{-k_i}E_i
\end{equation*}

The Henriques-Kamnitzer coboundary construction in \cite{MR2219257} uses the Lusztig
involution.
\begin{defn}
	The category $U_q(\sC)-mod$ is a coboundary category
	where $\sigma_{B,C}\colon B\otimes C \rightarrow C\otimes B$ is given by
	\begin{equation*}
	\sigma_{B,C}(b\otimes c) = \xi_{C\otimes B}(\xi_C(c)\otimes\xi_B(b)) 
	\end{equation*}
\end{defn}

These two constructions are shown to agree, up to signs, in \cite[Corollary 8.4]{MR2496310}.

\section{Crystals}
Crystals should be thought of as a discrete analogue of highest weight
representations. The most important property for this paper is a tensor
product rule. This tensor product rule is used to construct finite sets
which are discrete analogues of the isotypical components in tensor
products.

\subsection{Crystals}
For the cyclic sieving phenomenon we require a set with a bijection instead of
a vector space with an endomorphism. This is constructed using crystals.
Then we show that $N(\otimes^rV)$ has a basis which is permuted by the action
of $C_r$ and that this agrees with promotion on invariant words.

Here we give the definition of a crystals and the tensor product rule.
These were introduced in \cite{MR1090425}.
A good introduction to the theory of crystals, with examples, is \cite{MR1881971}.

\begin{defn} A \emph{crystal} is a set $B$ together with maps
 $e_\alpha,f_\alpha\colon B \rightarrow B\amalg \{0\}$ for $\alpha\in\Pi$
such that $e_\alpha b=b'$ if and only if $f_\alpha b'=b$.
\end{defn}
It is standard practice to represent a crystal, $B$, as a directed graph
with vertices $B$ and edges labelled by simple roots. If there is an edge
$b\rightarrow b'$ with $b,b'\in B$ and labelled by $\alpha$ then
$e_\alpha b = b'$ and $f_\alpha b' = b$.

A \emph{homomorphism} of crystals $\psi\colon B\rightarrow B'$ is a map
$\psi\colon B \coprod \{0\}\rightarrow B'\coprod \{0\}$ that commutes
with all $e_\alpha$ and all $f_\alpha$ and such that $\psi(0)=0$.

Define functions $\phi_\alpha,\varepsilon_\alpha\colon B\rightarrow\bZ$ by
\begin{equation*}
 \phi_\alpha(b)=\max_{k\ge 0}\{ k:f_\alpha^k(b)\ne 0\} \qquad
 \varepsilon_\alpha(b)=\max_{k\ge 0}\{ k:e_\alpha^k(b)\ne 0\}
\end{equation*}

Define the weight function $\wt\colon B\rightarrow P$ by
\begin{equation*}
 \langle \wt(b),\alpha\rangle = \phi_\alpha(b) - \varepsilon_\alpha(b)
\end{equation*}

Crystals have some basic properties which justify viewing them as a combinatorial
analogue of representations. First, a crystal corresponds to an irreducible representation if and only if the graph of the crystal is connected; second, the disjoint union of crystals corresponds to direct sum of representations and third, the analogue of Schur's lemma is that the set of homomorphisms between two connected crystals is a singleton if they are isomorphic and is empty otherwise.

For $\varpi\in P_+$, let $B(\varpi)$ be the connected crystal corresponding to $L_q(\varpi)$. These properties show that any crystal, $B$, has a canonical decomposition as
\begin{equation}
B \cong \coprod_{\varpi\in P_+} T(\varpi) \times B(\varpi)
\end{equation}
where $T(\varpi)$ is a set. These sets are discrete analogues of the isotypical subspaces and
the cardinality of the set is the dimension of the vector space.

Define the character of the crystal to be $\sum_{b\in B} \wt(b)$. This is an
element of the group ring $\bZ P$ and is invariant under the action of the Weyl
group. A representation and its crystal have the same character.

The tensor product rule for crystals is $B\otimes B'$ as a set is
$B\times B'$ and the raising and lowering operators are
\begin{align*}
 e_\alpha(b,b') &= \begin{cases}
(e_\alpha b,b') & \text{if $\varepsilon_\alpha(b) > \phi_\alpha(b')$} \\
(b,e_\alpha b') & \text{otherwise}
\end{cases} \\
 f_\alpha(b,b') &= \begin{cases}
(f_\alpha b,b') & \text{if $\varepsilon_\alpha(b) \ge \phi_\alpha(b')$} \\
(b,f_\alpha b') & \text{otherwise}
\end{cases}
\end{align*}

An important property of this tensor product is that it is associative.
The unit for the tensor product is the crystal with one element.

In particular, we can form the tensor powers of a crystal. For a crystal $B$,
the elements of $\otimes^rB$ are words of length $r$ in the elements of $B$.
This crystal is canonically decomposed as
\begin{equation}\label{eq:dec}
\otimes^r B \cong \coprod_{\varpi\in P_+}\Hom(B(\varpi),\otimes^r B) \times B(\varpi)
\end{equation}
This is a crystal isomorphism with each $\Hom(B(\varpi),\otimes^r B)$. This is a
discrete analogue of the decomposition of $\otimes^r M$ where $M$ is the representation
associated to $B$. It is common practice to identify $\Hom(B(\varpi),\otimes^r B)$
with the set of highest weight words in $\otimes^r B$ of weight $\varpi$. In
special cases these sets are identified with some version of tableaux and the
isomorphism \eqref{eq:dec} is described by an insertion algorithm.

In this paper we are mainly interested in the case $\varpi=0$. For $r>0$,
$\Hom(B(0),\otimes^r B)$ is the set of isolated vertices in the graph of $\otimes^r B$.
These are called the invariant words and this set is denoted $(\otimes^r B)_*$.

For a finite crystal $B$, the map $\xi\colon B\rightarrow B$ is determined by the
properties that $b$ and $\xi(b)$ are in the same connected component of $B$ and
\begin{equation*}
e_\alpha\xi(b)=\xi(f_{\theta(\alpha)}b)\qquad f_\alpha\xi(b)=\xi(e_{\theta(\alpha)}b)
\end{equation*}
and $\wt(\xi(b))=w_0(\wt(b))$.

\begin{defn}
	The category of finite crystals is a coboundary category
	where $\sigma_{B,C}\colon B\otimes C \rightarrow C\otimes B$ is given by
	\begin{equation*}
	\sigma_{B,C}(b\otimes c) = \xi_{C\otimes B}(\xi_C(c)\otimes\xi_B(b)) 
	\end{equation*}
\end{defn}

\subsection{Based modules}
Based modules were introduced in \cite[Chapter 27]{MR1227098} and are discussed in \cite{MR1637330}.

Define $\bA=\bQ[[q^{-1}]] \cap \bQ(q)$. Then we have a homomorphism 
$ev\colon \bA\rightarrow \bA / q^{-1}\bA \cong \bQ$. Loosely speaking, $\bA$ is the set of rational functions whose evaluation at $q^{-1}=0$ is defined and this evaluation is given by $ev$. Formally, $\bA$ is
a discrete valuation ring with maximal ideal generated by $q^{-1}$ and residue field $\bQ$.

\begin{defn} The data for a based module $(M,B)$ consists of $M\in U_q(\sC)-mod$, $B$ a basis of $M$ over $\bQ(q)$
and $\overline{\phantom{a}}$, a $\bQ$-linear involution of $M$.
\begin{description}
	\item[Weight spaces] The basis $B$ is required to be compatible with
the weight space decomposition of $M$, $M\cong \oplus_{\varpi\in P_+} M_\varpi$,  in the sense that 
\begin{equation}
B=\coprod_{\varpi\in P_+} (B\cap M_\varpi)
\end{equation}
\item[Bar involution] The involution is required to be compatible with the bar involution in the sense that
$\overline{ub}=\overline{u}b$ for $u\in U_q(\sC)$ and $b\in B$.
\item[Integral form] The basis $B$ gives an integral form in the sense that $\bZ[q,q^{-1}]B$ is
preserved by the integral form of $U_q(\sC)$.
\item[Crystal] The basis $B$ gives a crystal. Extend $ev$ to $ev\colon \bA B\rightarrow \bQ B$. The underlying set is $ev(B)$ and the maps
$e_\alpha,f_\alpha \colon ev(B)\rightarrow ev(B)\coprod \{0\}$ are defined by
$e_\alpha ev(b) = ev(\widetilde{E}_\alpha b)$ and $f_\alpha ev(b) = ev(\widetilde{F}_\alpha b)$
where the maps $\widetilde{E}_\alpha$ and $\widetilde{F}_\alpha$ are the Kashiwara operators.
\end{description}\end{defn}

A homomorphism of based modules $\psi\colon (M,B)\rightarrow (M',B')$ is a morphism in
$U_q(\sC)-mod$, $\psi\colon M\rightarrow M'$ such that $\psi(B)\subseteq B'\cup \{0\}$.

Based modules are closely related to canonical bases. The three original constructions of canonical bases in \cite{MR1090425}, \cite{MR1062796}, \cite{MR1035415} are all technical. More recently it is shown in
\cite{MR2630039} that the construction in \cite{MR2342692} gives the canonical bases. Every irreducible
representation with its canonical basis is a based module; also, every based module has a filtration such
that every subquotient is of this form.

There are some simple examples which can be constructed directly; the representations in Example \ref{ex:sl2},
the minuscule representations in \cite[Chapter 5A]{MR1359532} and the adjoint representations in \cite{adjoint}.

The basis $B$ is preserved by the Lusztig involution, \cite[Theorem 3.3]{MR1182165}. It follows that the mapping from based modules to crystals is a coboundary functor.

\subsection{Promotion}
Let $\varpi\in P_+$ and consider $L_q(\varpi)\in U_q(\sC)-mod$.
For $r>0$ there is an action of the cactus group $\cactus_r$ on $\otimes^r L_q(\varpi)\in U_q(\sC)-mod$. Promotion is given by the action of the element $\langle 1,r-1,r\rangle\in\cactus_r$, or equivalently,
by the map $\sigma_{B',B}$ with $B'=\otimes^{r-1}B$.
Promotion can then be considered as an operator on each isotypical subspace of $\otimes^r L_q(\varpi)$ and, in particular, as an operator on invariant tensors.

Drinfel$'$d promotion is the action of this element using the Drinfel$'$d coboundary structure and Henriques-Kamnitzer promotion is the action of this element using the Henriques-Kamnitzer coboundary structure.

Let $(\otimes^rB)_\ast$ be the set of isolated vertices for a crystal $B$. The Henriques-Kamnitzer coboundary structure defines a coboundary structure on crystals and hence promotion as a bijection on $(\otimes^rB)_\ast$.
This can also be described
without using the coboundary structure as follows: first, remove the first letter
(which is the highest weight letter) to leave a lowest weight word; second,
take the highest weight word in the same component (which is a copy of $B$); and
third, add the lowest weight letter to the end of the word.

For $\varpi\in P_+$, put $c(\varpi)=(\varpi,\varpi+2\rho)$ where $2\rho$ is the sum of the positive roots. Then the action of the quadratic Casimir on $L(\varpi)$ is the scalar operator $c(\varpi)$. The $q$-analogue is shown in Figure \ref{fig:twist}.
\begin{figure}
\begin{center}	
	\begin{tikzpicture}[line width=2pt, scale=0.75]
	\fill[color=blue!20] (-1.5,-1.5) rectangle (2,1.5);
	\draw (-1,1.5) -- (-1,1) -- (1,-1) arc(225:360:0.31) -- (1.5,0);
	\fill[color=blue!20] (0,0) circle(0.5);
    \draw (-1,-1.5) -- (-1,-1) -- (1,1) arc(135:0:0.31) -- (1.5,0);	
	\end{tikzpicture}
    \raisebox{1cm}{$\,=\pm q^{c(\varpi)}\,$}
	\begin{tikzpicture}[line width=2pt, scale=0.75]
	\fill[color=blue!20] (-0.5,-1.5) rectangle (0.5,1.5);
	\draw (0,-1.5) -- (0,1.5);
	\end{tikzpicture}	
\end{center}
\caption{Twist}\label{fig:twist}
\end{figure}
\begin{prop}\label{prop:b} Drinfel$'$d promotion on $N(\otimes^r L(\varpi))$
is given by the action of the long cycle.
\end{prop}

\begin{proof}
The braid group, $B_r$, acts on
$N(\otimes^r L_q(\varpi))$. Put $R=\sigma_1\sigma_2\dotsc \sigma_{r-1}\in B_r$
and then $R^{\mathrm{op}}=\sigma_{r-1}\dotsc \sigma_2 \sigma_{1}$.
Then Drinfel$'$d promotion is the unitarisation of $R$, as described in \S~\ref{sec:un}.

Let $\rho$ be rotation. Then, by a modification of Proposition \ref{lem:rot} using Figure \ref{fig:twist},
\begin{equation*}
R=\pm q^{c(\varpi)}\rho \qquad R^{\mathrm{op}}=\pm q^{c(\varpi)}\rho^{-1}\,.
\end{equation*}
where the two signs are the same. Therefore, $R^{\mathrm{op}}R=q^{2c(\varpi)}$,
$(R^{\mathrm{op}}R)^{-1/2}=q^{-c(\varpi)}$ and the unitarisation of $R$ is
$q^{-c(\varpi)}R=\pm\rho$. This is independent of $q$ and for $q=1$ is the action of the long cycle.
\end{proof}

\begin{prop}
Henriques-Kamnitzer promotion acting on $N(\otimes^r L_q(\varpi))$ is given by rotation.
\end{prop}

\begin{proof} The Drinfel$'$d coboundary structure and the Henriques-Kamnitzer coboundary structure
are shown to agree, up to sign, in \cite[Corollary 8.4]{MR2496310}. In particular, Drinfel$'$d promotion
and Henriques-Kamnitzer promotion agree, up to sign.
The action of the long cycle agrees with rotation, up to sign, by Proposition \ref{lem:rot}.
The two signs are the same so the result follows.
\end{proof}

\begin{prop}
Henriques-Kamnitzer promotion is represented on the canonical basis of $N(\otimes^r L_q(\varpi))$ by the permutation matrix of Henriques-Kamnitzer promotion acting on invariant words.
\end{prop}
\begin{proof} In \cite[\S~28.2.9]{MR1227098} it is shown that Henriques-Kamnitzer promotion is represented on the canonical basis of $N(\otimes^r L_q(\varpi))$ by some permutation matrix. Putting $q=\infty$ in this matrix gives the same matrix and this is the permutation matrix of Henriques-Kamnitzer promotion acting on invariant words.
\end{proof}

Then an immediate Corollary is:
\begin{thm}
Rotation is represented on the canonical basis of $N(\otimes^r L_q(\varpi))$ by the permutation matrix of Henriques-Kamnitzer promotion acting on invariant words.
\end{thm}
\subsection{Cyclic sieving phenomenon}
In this section we state our main theorem.

Let $\sC$ be a finite type Cartan matrix and $\varpi\in P_+$.
Associated to $\varpi$ is $L(\varpi)$, an irreducible representation, and $B(\varpi)$,
a connected crystal.
Then we construct a triple $(X,c,P)$ for each $r\ge 0$.
The set $X$ is the set of isolated vertices in $\otimes^rB(\varpi)$ and $c$ is promotion.
Put $\chi_r= \bch N(\otimes^r L(\varpi))$. Define $P$ to be $\fd \omega\chi_r$ if
$L(\varpi)$ has a symplectic form and to be $\fd \chi_r$ otherwise.

Proposition \ref{prop:b} and Theorem \ref{thm:fund} together imply our main theorem:
\begin{thm}\label{thm:main}
For all $\varpi$ and $r$, the triple $(X,c,P)$ exhibits the cyclic sieving phenomenon. 
\end{thm}

As an application, we have one of the main results of \cite{rhoades}.
The following triple exhibits the cyclic sieving phenomenon:
\begin{itemize}
 \item The set $X$ is the set of rectangular standard tableaux with $n$ rows and $k$ columns
 \item The bijection $c$ is promotion.
 \item The polynomial $P$ is $\fd (s_{n^k})$
\end{itemize}

In order to deduce this from Theorem \ref{thm:main} we take $M$ to be the vector
representation of $\fsl(n)$. The space of invariant tensors, $N(\otimes^rM)$, is
0 unless $n\mid r$ and for $r=kn$ the Frobenius character of $N(\otimes^rM)$ is
$s_{n^k}$. Let $B$ be the crystal of $M$.

The set $(\otimes^rB)_\ast$ is in
bijection with rectangular standard tableaux with $n$ rows and $k$ columns and this
bijection interwines the two promotion maps.

The examples given by taking $M$ to be a symmetric power of the vector representation
of a symplectic group are discussed in detail in \cite{3592}.

\section{Octonions}\label{sec:g2p}
The Lie algebra $G_2$ is a simple Lie algebra of dimension 14. The most straightforward
construction of this Lie algebra is that it is the derivation algebra of the octonions. 
There are two fundamental representations. One is the adjoint representation,
the other has dimension 7, which we refer to as the vector representation and which can
be taken to be the imaginary octonions.

The simple Lie algebra $G_2$ is the Lie algebra associated to the Cartan matrix
\begin{equation*}
\sC=\begin{pmatrix} 2 & -3 \\ -1 & 2 \end{pmatrix}
\end{equation*}
This Lie algebra has dimension 14 and the Weyl group is $D_6$, the dihedral group
of order 12.

\subsection{Invariant tensors}
In this section we give the results of computing the Frobenius characters of
$N(\otimes^rV)$ for $V$ one of the two fundamental representations and for $r$
small. There is no formula known for these Frobenius characters.

We record the multiplicity of the trivial representation in $\bS^\lambda(V)$ 
for $|\lambda|=2,3$ for these two representations in the following table
\begin{center}\begin{tabular}{r|cc|ccc}
& $[2]$ & $[1,1]$ & $[3]$ & $[2,1]$ & $[1,1,1]$ \\ \hline
vector & 1 & 0 & 0 & 0 & 1 \\
adjoint & 1 & 0 & 0 & 0 & 1 \\
\end{tabular}\end{center}
This table just records the information that each representation has an invariant	symmetric bilinear form and an invariant anti-symmetric trilinear form and there
are no other invariant tensors of degree at most $3$. For $r=2$ this gives the polynomial
$1$. For $r=3$ this gives the polynomial $q^3$. Then reducing modulo $q^3-1$ gives $1$.
	
For $r=4$ we have 
\begin{center}\begin{tabular}{r|ccccc}
& $[4]$ & $[3,1]$ & $[2,2]$ & $[2,1,1]$ &$[1,1,1,1]$ \\ \hline
vector & 1 & 0 & 1 & 0 & 1 \\
adjoint & 1 & 0 & 2 & 0 & 0 \\
\end{tabular}\end{center}
		
This gives the polynomials
\begin{center}\begin{tabular}{r|r|r}
vector &  $1+q^2+q^4+q^6$ & $2+2q^2$ \\
adjoint &  $1+2q^2+2q^4$ & $3+2q^2$
\end{tabular}\end{center}
This predicts that for the vector representation there are two orbits of size $2$
and for the adjoint representation there are two orbits of size $2$ and one orbit of size $1$.
			
For $r=5$ we have 
\begin{center}\begin{tabular}{r|ccccccc}
& $[5]$ & $[4,1]$ & $[3,2]$ & $[3,1^2]$ & $[2^2,1]$ & $[2,1^3]$ & $[1^5]$ \\ \hline
vector & 0 & 0 & 0 & 1 & 0 & 1 & 0 \\
adjoint & 0 & 0 & 0 & 2 & 0 & 1 & 0 \\
\end{tabular}\end{center}
				
This gives the polynomials
\begin{center}\begin{tabular}{r|r}
vector & $q^3+q^4+2q^5+2q^6+2q^7+q^8+q^9$ \\
adjoint &  $2q^3+2q^4+4q^5+3q^6+3q^7+q^8+q^9$
\end{tabular}\end{center}
and the reductions modulo $q^5-1$ are
\begin{center}\begin{tabular}{r|r}
vector &  $2+2q+2q^2+2q^3+2q^4$ \\
adjoint &  $4+3q+3q^2+3q^3+3q^4$
\end{tabular}\end{center}
This predicts that for the vector representation there are two free orbits 
and for the adjoint representation there are three free orbits and fixed point.
						
For $r=6$ we have 
\begin{center}\begin{tabular}{r|cccccc}
& $[6]$ & $[4,2]$ & $[3,2,1]$ & $[3,1^3]$ & $[2^3]$ & $[2,1^4]$  \\ \hline
vector & 1 & 1 & 0 & 1 & 2 & 1  \\
adjoint & 2 & 3 & 1 & 1 & 4 & 1  \\
\end{tabular}\end{center}
							
For the vector representation this gives the polynomial
\[ 1+q^2+q^3+2q^4+q^5+5q^6+2q^7+5q^8+4q^9+5q^{10}+2q^{11}+4q^{12}+q^{13}+q^{14} \]
For the adjoint representation this gives the polynomial
\[ 2+3q^2+3q^3+7q^4+5q^5+13q^6+7q^7+12q^8+8q^9+9q^{10}+3q^{11}+6q^{12}+q^{13}+q^{14} \]
and the reductions modulo $q^6-1$ are
\begin{center}\begin{tabular}{r|r}
vector &  $10+3q+7q^2+5q^3+7q^4+3q^5$ \\
adjoint &  $21+8q+16q^2+11q^3+16q^4+8q^5$
\end{tabular}\end{center}
For the vector representation this corresponds to 3,4,2,1 orbits of sizes 6,3,2,1 respectively. For the adjoint representation this corresponds to 8,8,3,2 orbits of sizes 6,3,2,1 respectively.

\subsection{Webs}\label{sec:web}
In this section we discuss webs. These have only been constructed in a few
cases. These diagrams can be rotated and we conjecture that this agrees with
the rotation map. For $SL(2)$ these are the diagrams in \cite{MR1446615}
and this is an elementary construction of the based modules and their tensor
products. The web diagrams for rank two simple Lie algebras are constructed in
\cite{MR1403861}. The cyclic sieving phenomenon in \cite{rhoades} in the cases
$A_1=\SL(2)$ and $A_2=\SL(3)$ are described in terms of diagrams in \cite{MR2519848}.
This was the inspiration for this paper.

In this section we confirm the predictions of \S \ref{sec:oct} using the web bases
of \cite{MR1403861}.

A \emph{trivalent planar graph} means a subset $\Gamma$ of the lower half plane $H$ such that
every point of $\Gamma$ has an open neighbourhood in $H$ which is homeomorphic to
one of the following three cases:

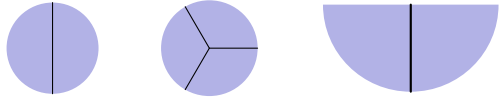
\[
\includegraphics[height=.5in]{csp.1}\qquad
\includegraphics[height=.5in]{csp.2}\qquad
\includegraphics[height=.5in]{csp.3} \]

The centre point in the second case is called a trivalent vertex.
In the third case the straight boundary is in the boundary of $H$.
Each of these points is a boundary point.

A trivalent planar graph is \emph{non-positive} if it does not contain any
of the following forbidden configurations

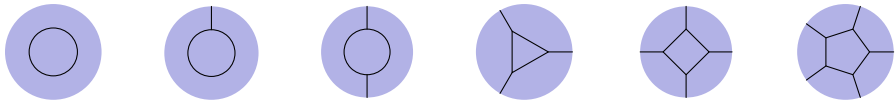
\[
\includegraphics[height=.5in]{csp.4}\qquad
\includegraphics[height=.5in]{csp.5}\qquad
\includegraphics[height=.5in]{csp.6}\qquad
\includegraphics[height=.5in]{csp.7}\qquad
\includegraphics[height=.5in]{csp.8}\qquad
\includegraphics[height=.5in]{csp.9}
\]

The set $X(r)$ is the set of non-positive trivalent planar graphs with $r$ boundary
points. Note that there is no restriction
on the number of trivalent vertices of the graph or on the number of connected
components.
A curvature argument using the isoperimetric inequality for surfaces of non-positive
curvature shows that for all $r\ge 0$, the set $X(r)$ is finite. The numbers
$|X(r)|$ for $0\le r\le 9$ are given in the following table
\begin{center}\begin{tabular}{c|rrrrrrrrr}
		$r$ & 0 & 1 & 2 & 3 & 4 & 5 & 6 & 7 & 8  \\
		$|X(r)|$ & 1 & 0 & 1 & 4 & 10 & 35 & 120 & 455 & 1792
	\end{tabular}\end{center}
This is sequence \href{http://www.research.att.com/~njas/sequences/A059710}{A059710} in \cite{OEIS}.
	
The diagrams in this section are drawn with two types of edge;
namely a single edge and a double edge. The conversion to non-positive trivalent planar graphs
is given by making the following substitution for each double edge:
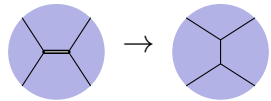
\[ \incg{csp.19} \rightarrow \incg{csp.20} \]
	
The cyclic group of order $r$ acts on $X(r)$ by rotation. 
For $r=0$ the only diagram is the empty diagram. There is no diagram
for $r=1$. For $r=2,3$ there is one diagram. These are 
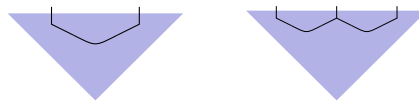
\[ \includegraphics[height=.5in]{triangle.10}\qquad 
\includegraphics[height=.5in]{triangle.11}\]
\begin{ex}
For $r=4$ we have four diagrams. The two orbits of order two are:
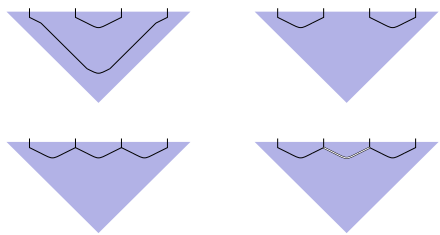
\[ \includegraphics[height=.5in]{triangle.12}\qquad 
\includegraphics[height=.5in]{triangle.13}\]
\[ \includegraphics[height=.5in]{triangle.14}\qquad 
\includegraphics[height=.5in]{triangle.15}\]
\end{ex}
\begin{ex}
For $r=5$ we have ten diagrams. There are two orbits of order five. One is the orbit
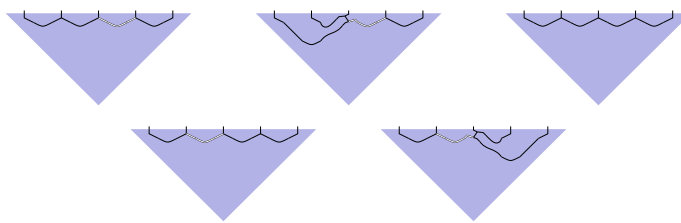
\[ \includegraphics[height=.5in]{triangle.16}\qquad 
\includegraphics[height=.5in]{triangle.17}\qquad
\includegraphics[height=.5in]{triangle.20}\]
\[ \includegraphics[height=.5in]{triangle.22}\qquad
\includegraphics[height=.5in]{triangle.25}\]	The other is the orbit
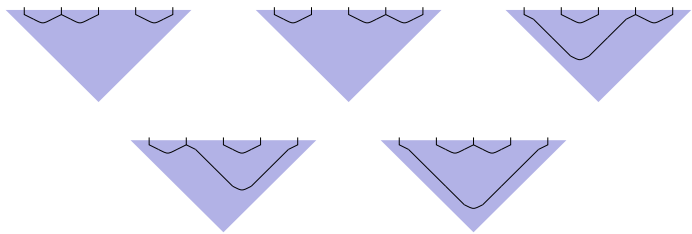
\[\includegraphics[height=.5in]{triangle.18}\qquad
\includegraphics[height=.5in]{triangle.19}\qquad
\includegraphics[height=.5in]{triangle.21}\]
\[ \includegraphics[height=.5in]{triangle.23}\qquad 
\includegraphics[height=.5in]{triangle.24}\]
\end{ex}
	
\subsection{Promotion}\label{sec:oct}

By iterating the tensor product we can construct the tensor powers of a crystal.
The vertices of $\otimes^rB$ are words of length $r$ in the alphabet of vertices
of $B$. The standard way to represent a word is as a sequence 
$(b_1,b_2,\dotsc ,b_r)$. An alternative is to represent a word as a set partition
of $\{1,2,\dotsc ,r\}$ into blocks labelled by $B$. In this representation each
vertex of the crystal is labelled by a subset of $\{1,2,\dotsc ,r\}$
and the effect of applying the raising operator $e_i$ to that a number moves along
an edge labelled $i$.

In this section we give examples of promotion for the seven dimensional
representation of the exceptional simple Lie algebra $G_2$.
In \cite{MR2320368} we constructed a flow diagram for each word.
In these examples we give three representation of each word;
the standard representation, the flow diagram and the labelling of the vertices
of the crystals by subsets of $\{1,2,\dotsc ,r\}$.

The crystal has seven vertices which we will label by the seven letters, $\{A,B,C,o,c,b,a\}$.
Words will then mean words in this alphabet.

The crystal is then given by
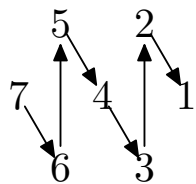
\[ \includegraphics[width=1in]{csp.21} \]
where vertices are positioned in the weight lattice.

The conventional way of drawing this crystal is
\begin{equation}
A\overset{\alpha}{\rightarrow}
B\overset{\beta}{\rightarrow}
C\overset{\alpha}{\rightarrow}
o\overset{\alpha}{\rightarrow}
c\overset{\beta}{\rightarrow}
b\overset{\alpha}{\rightarrow}
a
\end{equation}
where $\alpha$ is the short root and $\beta$ is the long root.

\begin{figure}
\[ \includegraphics[width=3in]{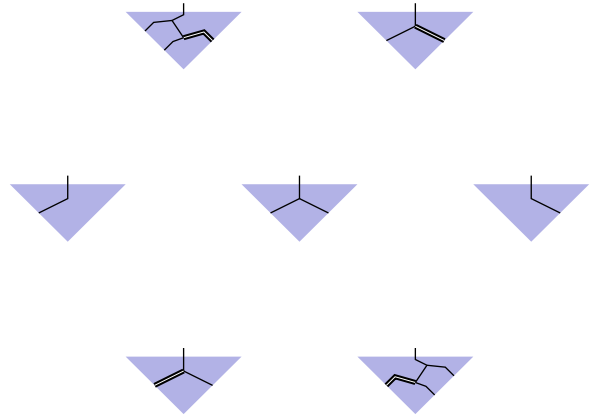} \]
\caption{Triangular diagrams}\label{fig:tr}.
\end{figure}

We will also represent words by the triangular diagrams introduced in \cite{MR2320368}. In this representation the seven vertices are labelled as in Figure \ref{fig:tr}.

The third way we will represent words of length $r$ is by a set partition (with empty blocks allowed) of $\{1,2,\dotsc ,r\}$ with blocks indexed by the seven vertices. In this representation each
vertex of the crystal is labelled by a subset of $\{1,2,\dotsc ,r\}$
and the effect of applying the raising operator $e_\alpha$ to that a number moves along
an edge labelled $\alpha$.

\begin{ex} We take the word aoA whose planar trivalent graph has a single
trivalent vertex and so is invariant under rotation. We remove the initial a to get
oA. Then we apply the lowering operators to get the following sequence of words
\[ oA\rightarrow cA\rightarrow bA\rightarrow bB\rightarrow aB\rightarrow aC\rightarrow ao \]
Then we add the A at the end to get aoA. Alternatively, the word aoA is represented by
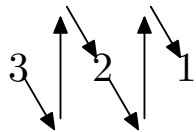
\[ \includegraphics[width=1in]{csp.26} \]
Then the sequence of words is represented by the sequence
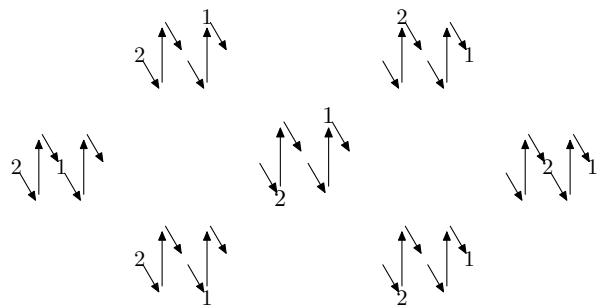
\[ \begin{array}{ccccc}
 & \includegraphics[width=.5in]{csp.29} & &
 \includegraphics[width=.5in]{csp.32} & \\
 \includegraphics[width=.5in]{csp.27} & &
 \includegraphics[width=.5in]{csp.30} & &
 \includegraphics[width=.5in]{csp.33} \\
 & \includegraphics[width=.5in]{csp.28} & &
 \includegraphics[width=.5in]{csp.31} &
\end{array} \]
\end{ex}

\begin{ex} We start with the word abcAoA. We remove the initial a to get
bcAoA. Then we apply the lowering operators to get the following vertices 
of a crystal
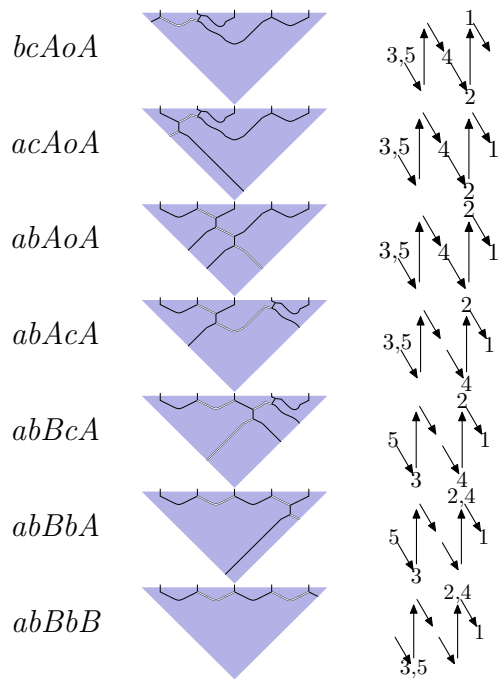
\begin{center}\begin{tabular}{ccc}
bcAoA & \incg{triangle.65} & \incg{csp.34} \\
acAoA & \incg{triangle.66} & \incg{csp.35} \\
abAoA & \incg{triangle.67} & \incg{csp.36} \\
abAcA & \incg{triangle.68} & \incg{csp.37} \\
abBcA & \incg{triangle.69} & \incg{csp.38} \\
abBbA & \incg{triangle.70} & \incg{csp.39} \\
abBbB & \incg{triangle.71} & \incg{csp.40} 
\end{tabular}\end{center}
Then we add the A at the end to get abBbBA. This gives $abcAoA \mapsto abBbBA$.
\end{ex}

\begin{ex} We start with the word abBbBA. We remove the initial a to get
bBbBA. Then we apply the lowering operators to get the following vertices 
of a crystal
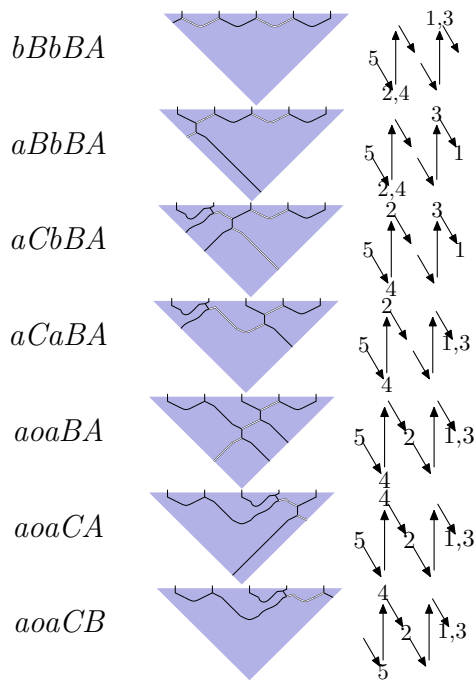
\begin{center}\begin{tabular}{ccc}
bBbBA & \incg{triangle.72} \incg{csp.41} \\
aBbBA & \incg{triangle.73} \incg{csp.42} \\
aCbBA & \incg{triangle.74} \incg{csp.43} \\
aCaBA & \incg{triangle.75} \incg{csp.44} \\
aoaBA & \incg{triangle.76} \incg{csp.45} \\
aoaCA & \incg{triangle.77} \incg{csp.46} \\
aoaCB & \incg{triangle.78} \incg{csp.47}
\end{tabular}\end{center}
Then we add the A at the end to get aoaCBA. This gives
$abBbBA\mapsto aoaCBA$.
\end{ex}

These examples support the following conjecture:
\begin{conj} The bijection in \cite{MR2320368} between non-positive trivalent
	graphs with $r$ boundary points and isolated words in $\otimes^rB$ interwines
	rotation and promotion.
\end{conj}

\section{Energy statistic}
A triple $(X,c,P)$ which exhibits the cyclic sieving phenomenon is more interesting from a combinatorial perspective if $P$ is the generating function of a statistic on $X$. That is, there is a function $\mathrm{st}\colon X\rightarrow \bN$ such that $P=\sum_{x\in X} q^{\mathrm{st}(x)}$.

In this section we show that for some representations $P$ is a classically restricted one dimensional configuration sum. These polynomials are given by the fermionic formula and are the generating function for a statistic known as energy. The representations for which this holds are listed in section \ref{sec:cir}. This construction of examples of the cyclic sieving phenomenon still uses Proposition \ref{prop:b} but avoids Theorem \ref{thm:fund}.

The energy function is purely combinatorial but the representation theory background is the study
of Kirillov-Reshetikhin modules. This is a large and difficult subject with its origins in the
study of exactly solvable models in statistical mechanics.

\subsection{Kirillov-Reshetikhin modules}
Let $\fg$ be a semisimple Lie algebra. Then let $\widetilde{\fg}$ be the affine algebra
of $\fg$ and $\widetilde{\fg}'$ be its derived subalgebra. The affine algebra is infinite
dimensional and the derived algebra has codimension one. The affine algebra is a Kac-Moody
algebra and the derived algebra is defined by omitting the grading operator. The affine
algebra has no non-trivial finite dimensional representations, whereas the derived algebra has an
interesting category of finite dimensional representations. These both have quantised enveloping
algebras. Our interest is in $U_q(\widetilde{\fg}')$ and its category of finite dimensional representations.

The Yangian, $Y(\fg)$, is a deformation of $U(\fg[t])$ associated with a Lie bialgebra structure
on $\fg[t]$. The category of finite dimensional dimensional representations of $U_q(\widetilde{\fg}')$
should be understood as a $q$-analogue of the category of finite dimensional dimensional representations of $Y(\fg)$ even though it is not obvious that it is given by a deformation quantisation. One way this
is made manifest is that a representation of $Y(\fg)$ gives a rational solution to the Yang-Baxter
equation (with spectral parameter) and a representation of $U_q(\widetilde{\fg}')$ gives a trigonometric solution. 

The most important, and best understood, representations of these quantum groups are
the Kirillov-Reshetikhin modules (aka KR-modules) and their tensor products.
There is a KR-module for each simple root $\alpha$ and each $s>0$. The associated module
for the Yangian is denoted $W^{(\alpha,s)}$ and for the quantum affine algebra by
$W_q^{(\alpha,s)}$.
We refer the reader to the excellent survey, \cite{MR2642561}, for more information.

The property of the Kirillov-Reshetikhin modules we are interested in is that they have a
crystal. In fact, these modules and their tensor products are the only known examples
of finite dimensional representations of the quantum affine algebra which admit a crystal.
The crystal of $W_q^{(\alpha,s)}$ is denoted $B^{(\alpha,s)}$. Two uniform
constructions of these crystals are the alcove path model \cite{MR3090977} and the 
level zero projection of Littelmann paths \cite{MR1356780}.

\subsection{Classically irreducible modules}\label{sec:cir}
The pair $(\alpha,s)$ is classically irreducible if the classical restriction
of $B^{(\alpha,s)}$ is connected. Equivalently, the restriction of $W_q^{(\alpha,s)}$
to $U_q(\fg)$ is irreducible. In this case the restriction is the
highest weight representation $V(s\omega_\alpha)$. 

First we list the fundamental representations $V(\omega_\alpha)$ such that $W_q^{(\alpha,s)}$ is classically
irreducible for all $s\ge 1$. For type $A$ these are all of the fundamental representations.
For the remaining types the list is
\begin{description}
	\item[Type $B$] the defining representation
	\item[Type $C$] the extreme fundamental representation whose dimension is a Catalan number
	\item[Type $D$] the defining representation and both half-spin representations
	\item[Type $E_6$] the dual pair of twenty seven dimensional representations
	\item[Type $E_7$] the fifty six dimensional representation
\end{description}
These are all perfect. These are also the nodes of the Dynkin diagram such that
the corresponding coordinate of the null root, $\delta$, is $1$.

Next we list the fundamental representations $V(\omega_\alpha)$ such that $W_q^{(\alpha,1)}$ is classically
irreducible but $W_q^{(\alpha,2)}$ is classically reducible.
\begin{description}
	\item[Type $B$] the spin representation
	\item[Type $C$] the defining representation
	\item[Type $G_2$] the seven dimensional representation
	\item[Type $F_4$] the twenty six dimensional representation
\end{description}
None of these are perfect.

Note that $E_8^{(1)}$ has no classically irreducible representation.

\subsection{Current algebras}
Let $\fg$ be a simple Lie algebra. The current algebra is the Lie algebra
$\fg[t] = \fg\otimes \bC[t]$. The Lie bracket is
\begin{equation*}
[x\otimes f(t),y\otimes g(t) ] = [x,y]\otimes f(t)g(t)
\end{equation*}
This has a grading where $t$ has degree 1 and this gives a grading on the universal
enveloping algebra $U(\fg[t])$. Let $V$ be a finite dimensional representation of
$U(\fg[t])$ together with a cyclic vector $v\in V$. This gives a filtration on $V$ by
\begin{equation*}
V^{(i)} = \sum_{j\le i} U(\fg[t])^j .v
\end{equation*}

The Yangian and the current algebra both have a family of automorphisms parametrised by $\bC$.
For the current algebra, the automorphisms $\tau_c$, for $c\in\bC$, are given by
\begin{equation*}
\tau_c\colon x\otimes f(t) \mapsto x\otimes f(t-c)
\end{equation*}
Then each representation, $V$, gives a family of representations by $V(c) = \tau_c^* V$.

Now let $V_k$ for $1\le k\le r$ be a sequence of representations each with a
vector $v_k$ which generates the module. Then the tensor product $V_1(u_1)\otimes \dotsb \otimes V_r(u_r)$ has the cyclic
vector $v_1\otimes \dotsb \otimes v_r$ provided the parameters are distinct.
The fusion product, $V_1(u_1)\ast \dotsb\ast V_r(u_r)$, is the graded representation
associated to the filtration on the tensor product representation with this choice of
cyclic vector.

For any representation, $V$, of $\fg$ we have a family of representations of $\fg[t]$,
$V(c)$ for $c\in\bC$, where
\begin{equation*}
(x\times f(t)).v = f(c) x.v
\end{equation*}

For any sequence $\brs$ and any distinct values of the parameters consider
$W^{(\alpha_1,s_1)}(u_1)\ast \dotsb \ast W^{(\alpha_r,s_r)}(u_r)$ as a graded 
$\fg$-module. Then each isotypical component is a graded vector space. For a graded
vector space $V=\oplus_k V^{[k]}$ the graded dimension is the polynomial
\begin{equation*}
\sum_k \dim (V^{[k]}) q^k
\end{equation*}
Then for all $\lambda\in P_+$ and all sequences $\brs$ define $C(\lambda,\brs)$
to be the graded dimension of the isotypical component 
\begin{equation*}
\Hom_\fg( V(\mu) , W^{(\alpha_1,s_1)}(u_1)\ast \dotsb \ast W^{(\alpha_r,s_r)}(u_r)
\end{equation*}

It is shown in \cite{MR1729359} that for all $\lambda\in P_+$ and all sequences $\brs$ the polynomial $C(\lambda,\brs)$ is independent of the choice of $u_1,u_2,\dots ,u_r$.

Let $\zeta$ be a primitive $r$-th root of unity.
\begin{prop}\label{prop:lc} The long cycle acts by $\zeta^k$ on the component of degree $k$ of
the fusion product $V(1)\ast V(\zeta) \ast \dotsb \ast V({\zeta^{r-1}})$.
\end{prop}
The proof is from \cite[\S 4.2]{MR3177928}.
\begin{proof}
Let $\rho$ be the action of the long cycle given by
\begin{equation*}
\rho(v_1\otimes \dotsb \otimes v_r) = v_2 \otimes \dotsb \otimes v_r \otimes v_1
\end{equation*}
Then $\rho$ commutes with the action of $\fg$ but not with the action of $\fg[t]$.
Instead we have
\begin{equation*}
x\otimes t^k \circ\rho = \zeta^k \rho\circ x\otimes t^k
\end{equation*}
\end{proof}

Then combining Proposition \ref{prop:b} and Proposition \ref{prop:lc} gives
the following example of the cyclic sieving phenomenon.
\begin{theorem}\label{thm:en} Let $W(\alpha,s)$ be classically irreducible.
Put $X=(\otimes^rB(\omega_\alpha))_\ast$ with $c$ acting by promotion.
Then $(X,c,P)$ exhibits the cyclic sieving phenomenon with $P=C(0,\brs)$ where
$\brs$ is the sequence with $(\alpha_i,s_i)=(\alpha,s)$ for $1\le i\le r$.
\end{theorem}

\subsection{Energy function}
Theorem \ref{thm:en} is unsatisfactory as it does not give any method for computing
the polynomial. The polynomials $C(\lambda,\brs)$ have several other interpretations.
One of the interpretations is as the generating function for the energy function.
As a polynomial this is given by the fermionic formula and this is an efficient
formula for computing these polynomials.

There are two ways of constructing the energy function on highest weight words.
Each involves some intricate
combinatorics. One way is to define the energy function as a statistic on rigged
configurations. It is a simple observation that the generating function is given
by the fermionic formula. Indeed, this was the motivation for introducing rigged
configurations. Then, in order to view the energy function as a statistic on
highest weight words, a bijection between rigged configurations and highest weight
words is required. The alternative is to define the energy function as a statistic
on highest weight words directly. This approach uses the combinatorial $R$-matrices.

\begin{ex} The polynomials in section \ref{sec:g2p} can be computed in sage using
\begin{verbatim}
sage: RC = RiggedConfigurations(['G', 2, 1], [[1,1]]*4 )
sage: mg = RC.module_generators
sage: [ a.cc() for a in mg if a.weight() == 0 ]
[8, 6, 6, 4]
\end{verbatim}	
\end{ex}
The result that the polynomials $C(\lambda,\brs)$ are given by the fermionic formula
has been proved not as an isolated result but in conjunction with several other
conjectures.
These include the conjectures set out in the remarkable paper \cite{MR1903978}.
We refer the interested reader to \cite{MR2767945} which gives an
overview.

\printbibliography

\end{document}